\documentclass[12pt,leqno,tbtags]{amsart}

\usepackage[T1]{fontenc}
\usepackage{amssymb}

\theoremstyle{plain}
\newtheorem{theorem}{Theorem}[section]
\newtheorem{proposition}[theorem]{Proposition}
\newtheorem{lemma}[theorem]{Lemma}
\newtheorem{corollary}[theorem]{Corollary}

\theoremstyle{definition}
\newtheorem{definition}[theorem]{Definition}

\theoremstyle{remark}
\newtheorem{remark}[theorem]{Remark}

\numberwithin{equation}{section}

\newcommand{\Cyclic}{\mathop{\text{\Large$\mathfrak S$}\vrule width 0pt
depth 2pt}}

\newcommand{\Lie}[1]{\textsl{#1}}
\newcommand{\lie}[1]{\mathfrak{#1}}

\DeclareMathOperator{\SO}{\Lie{SO}}
\DeclareMathOperator{\so}{\lie{so}}
\DeclareMathOperator{\SP}{\Lie{Sp}}
\DeclareMathOperator{\sP}{\lie{sp}}
\DeclareMathOperator{\SU}{\Lie{SU}}

\newcommand{\hol}{\lie{hol}}

\newcommand{\QK}[1]{\mathcal{QK}_{#1}}
\newcommand{\LC}{\nabla}
\newcommand{\Nt}{\widetilde\nabla}
\newcommand{\ow}[3]{#1{\otimes}#2{\wedge}#3}
\newcommand{\A}{\alpha}
\newcommand{\B}{\beta}
\newcommand{\Cc}{\gamma}

\newcommand{\JJ}{\mathfrak J}
\newcommand{\vo}[3]{#1{\vee}#2{\otimes}#3}

\newcommand{\qrad}{\varrho}
\newcommand{\Rt}{\widetilde R}

\newcommand{\V}{\mathcal V}
\newcommand{\VH}{\check{\mathcal V}}
\newcommand{\VE}{\hat{\mathcal V}}

\DeclareMathOperator{\Ad}{Ad}
\DeclareMathOperator{\ad}{ad}
\DeclareMathOperator{\diag}{diag}
\DeclareMathOperator{\Span}{\texttt{Span}}
\DeclareMathOperator{\Id}{Id}
\DeclareMathOperator{\re}{Re}
\DeclareMathOperator{\im}{Im}
\DeclareMathOperator{\HH}{\mathbb HH}
\DeclareMathOperator{\HP}{\mathbb HP}

\makeatletter
\newcommand{\inp}{\@ifstar{\inpFixed}{\inpScaled}}
\newcommand{\inpFixed}[2]{\langle #1, #2 \rangle}
\newcommand{\inpScaled}[3][]{\def\@tmp{#1}\def\@tmpd{}%
\ifx\@tmp\@tmpd\def\@tmpl{left}\def\@tmpr{right}%
\else\def\@tmpl{#1l}\def\@tmpr{#1r}\fi
\csname \@tmpl \endcsname \langle #2, #3 \csname \@tmpr \endcsname \rangle}
\newcommand{\norm}{\@ifstar{\normFixed}{\normScaled}}
\newcommand{\normFixed}[1]{\lVert #1 \rVert}
\newcommand{\normScaled}[2][]{\def\@tmp{#1}\def\@tmpd{}%
\ifx\@tmp\@tmpd\def\@tmpl{left}\def\@tmpr{right}%
\else\def\@tmpl{#1l}\def\@tmpr{#1r}\fi
\csname \@tmpl \endcsname \lVert #2 \csname \@tmpr \endcsname \rVert}
\newcommand{\abs}{\@ifstar{\absFixed}{\absScaled}}
\newcommand{\absFixed}[1]{\lvert #1 \rvert}
\newcommand{\absScaled}[2][]{\def\@tmp{#1}\def\@tmpd{}%
\ifx\@tmp\@tmpd\def\@tmpl{left}\def\@tmpr{right}%
\else\def\@tmpl{#1l}\def\@tmpr{#1r}\fi
\csname \@tmpl \endcsname \lvert #2 \csname \@tmpr \endcsname \rvert}
\makeatother

\begin{document}

\title[Homogeneous quaternionic K\"ahler structures]{Homogeneous quaternionic
K\"ahler structures and quaternionic\\ hyperbolic space}

\author{M. Castrill\'on L\'opez}
\address[MCL]{Departamento de Geometr\'\i a y Topolog\'\i a\\
Facultad de Matem\'aticas\\
Av. Complutense s/n\\
28040 Madrid\\
Spain} \email{mcastri@mat.ucm.es}

\author{P. M. Gadea}
\address[PMG]{Institute of Mathematics and Fundamental Physics\\
CSIC\\
Serrano 123\\
28006 Madrid\\
Spain} \email{pmgadea@iec.csic.es}

\author{A. F. Swann}
\address[AFS]{Department of Mathematics \& Computer Science\\
University of Southern Denmark\\
Campusvej 55\\
DK-5230 Odense M\\
Denmark} \email{swann@imada.sdu.dk}

\date{}


\thanks{Partially supported by DGICYT, Spain, under Grant BFM 2002-00141.
AFS partially supported by the EDGE, Research Training Network
HPRN-CT-2000-0010, of The European Human Potential Programme.}

\begin{abstract}
  An explicit classification of homogeneous quaternionic K\"ahler
  structures by real tensors is derived and we relate this to the
  representation-theoretic description found by Fino.  We then show how
  the quaternionic hyperbolic space $\HH(n)$ is characterised by admitting
  homogeneous structures of a particularly simple type.  In the process we
  study the properties of different homogeneous models for $\HH(n)$.
\end{abstract}

\maketitle

\begin{center}
  \begin{minipage}{0.8\textwidth}
    \begin{tiny}
      \tableofcontents
    \end{tiny}
  \end{minipage}
\end{center}
\newpage

\section{Introduction}
Representation theory has been successfully applied to the classification
of various geometric structures on differentiable manifolds in a number of
different settings, inspired by the initial work of Gray \& Hervella
\cite{Gray-H:16} for almost-Hermitian structures.

In the present paper, we first give a classification of homogeneous
quaternionic K\"ahler structures by real tensors.  Quaternionic K\"ahler
manifolds are characterised by having holonomy in $\SP(n)\SP(1)$,
$n\geqslant2$, and are Riemannian manifolds whose metrics are Einstein.
Many homogeneous examples are known following the work of
Wolf~\cite{Wolf:quaternionic}, Alekseevsky~\cite{Alekseevsky:solvable}, de
Wit \& van Proeyen \cite{DeWit-vP:special} and
Cort\'es~\cite{Cortes:Alekseevskian}, although a full classification has
not yet be found.  We hope that some of the techniques of this paper will
eventually lead to progress on this problem. A brief summary of the r\^ole
played by homogeneous quaternionic K\"ahler structures in theoretical
physics is provided at the end of this introduction.

Homogeneous Riemannian structures were studied systematically by Ambrose \&
Singer \cite{Ambrose-S:homogeneous} and Tricerri \& Vanhecke
\cite{Tricerri-V:homogeneous} in terms of tensors on manifolds.
In~\cite{Fino:torsion}, Fino specialised their results to the quaternionic
K\"ahler case, giving the abstract representation-theoretic decomposition
of the space~$\V$ of tensors satisfying the same symmetries as a
homogeneous quaternionic K\"ahler structure.  Our first main result is a
concrete description of this decomposition in terms of real tensors.  After
establishing some preliminaries in~\S\ref{sec:pre}, we give
(\S\ref{sec:cla}) a concrete orthogonal decomposition of~$\V$ into five
subspaces $\QK1$, \dots, $\QK5$ invariant under the action of
$\SP(n)\SP(1)$ and which we then relate to Fino's in Theorem~\ref{tchqKs}.

The first three modules in~$\V$ are distinguished by the fact that their
dimensions depend only linearly on $n=\dim M/4$.  We therefore single these
out for special attention.  By studying the interaction of the homogeneous
tensor with the quaternionic curvature we find in \S\ref{sec:geo} that
non-trivial homogeneous quaternionic K\"ahler structures in $\QK{1+2+3}$
are necessarily of type $\QK3$ and that a manifold admitting such a
structure has the curvature of quaternionic hyperbolic space.

This prompts us to study different homogeneous models for quaternionic
hyperbolic space~$\HH(n)$ in~\S\ref{sec:typ}.  We first recall how Witte's
refined Langlands decomposition of parabolic subgroups may be used to
determine all the connected groups acting transitively on a non-compact
symmetric space, and specialise to the case of~$\HH(n)$.  From the point of
view of Lie groups, the simplest model of $\HH(n)$ is as the solvable group
$AN$ in the Iwasawa decomposition $\SP(n,1)=KAN$.  However, the tensorial
description of this structure turns out to be rather complicated, being of
type $\QK{1+3+4}$.  The trivial homogeneous tensor corresponds to the
description of $\HH(n)$ as the Riemannian symmetric space
$\SP(n,1)/(\SP(n)\times\SP(1))$.  We find that the structures of type
$\QK3$ arise from a particular homogeneous description of~$\HH(n)$ as
$\SP(1)AN/\SP(1)$ with the isotropy representation depending on a positive
real parameter~$\lambda$.  In addition to the Lie-theoretic approach, we
provide a concrete description of this geometry on the open unit ball
in~$\mathbb H^n$.  The results contrast strongly with the case of real
hyperbolic space studied by Tricerri \& Vanhecke
\cite{Tricerri-V:homogeneous}, where the description as a solvable group is
particularly simple.

Combining the computations and constructions of sections \S\ref{sec:geo}
and \S\ref{sec:typ} we arrive at the following characterisation
of~$\HH(n)$.

\begin{theorem}
  \label{thm:main}
  A connected, simply-connected, complete quaternionic K\"ahler manifold of
  dimension $4n \geqslant 8$ admits a non-vanishing homogeneous
  quaternionic K\"ahler structure in the class $\QK{1+2+3}$ if and only if
  it is the quaternionic hyperbolic space.

  In this case, the homogeneous structure is necessarily of type~$\QK3$.
\end{theorem}

Earlier versions of some of the results of this paper were announced
in~\cite{Castrillon-Lopez--GS:qK-linear}.

We recall that hyperK\"ahler and quaternionic K\"ahler spaces appear in
various contexts in field and string theory.  For instance, they are found
in the formulation of the coupling of matter fields in $N = 2$
supergravity; that is, couplings of $n$ spin multiplets to supergravity
with two independent supersymmetric transformations.  Each multiplet
consists of $4$ real scalars and $2$ Majorana spinor fields.  The $4n$ real
scalars parameterise a $4n$-dimensional real manifold $M$ endowed with a
Riemannian metric such that the kinetic part of the Lagrangian reads as a
non-linear sigma model from the space time to $M$; i.e., a harmonic
Lagrangian.  This manifold is called the target manifold of the model.

Physical and topological considerations force the holonomy group of $M$ to
be a subgroup of~$\SP(n)$ (that is, $M$ is hyperK\"ahler) if the gravity
is considered as a background field, or $\SP(n)Sp(1)$ (that is, $M$ is
quaternionic K\"ahler) if the gravity is considered as a dynamical field.
The former case is called global supersymmetry and the latter local
supersymmetry.  We refer the reader to, for example,
\cite{Alekseevsky-C:isometry,Bagger-W:couplings,Cecotti:T,DeWit-vP:special}.

On the other hand, the only sigma models known to be integrable are those
whose target manifold is a homogeneous manifold (for example,
see~\cite{Bordemann-FLS:Lie-Poisson}).  Therefore, it seems reasonable the
existence of links between the classification of homogeneous structures and
a possible classification of some physical structures and models.  In fact,
it would be of interest to translate into physical terms each of the
classes $\QK i$ (cf.\ Theorem \ref{tchqKs}), the classes obtained by direct
sum of these and the corresponding mathematical structures.  Moreover, as
$N = 2$ supersymmetric non-linear sigma models require non-compact target
manifolds, the characterisation of quaternionic K\"ahler hyperbolic space,
a paradigm of non-compact spaces, in terms of homogeneous structures could
reveal interesting information of the aforementioned translation.

\section{Preliminaries}
\label{sec:pre}

\subsection{Ambrose-Singer equations}
\label{sec:AS}

Let $(M,g)$ be a connected, simply-connected, complete Riemannian manifold.
Ambrose \& Singer \cite{Ambrose-S:homogeneous} gave a characterisation for
$(M,g)$ to be homogeneous in terms of a $(1,2)$ tensor field~$S$.  The
tensor~$S$ is usually called a homogeneous Riemannian structure, and a
thorough study of these was made by Tricerri \&
Vanhecke in \cite{Tricerri-V:homogeneous} and a series of papers
by these authors and their collaborators.  If $\LC$ denotes the Levi-Civita
connection and $R$ its curvature tensor, then one introduces the torsion
connection $\Nt =\LC-S$ which satisfies the Ambrose-Singer equations
\begin{equation}
  \label{eq:AS}
  \Nt g  =  0, \quad \Nt R  =  0, \quad \Nt S  = 0.
\end{equation}

The manifold~$(M,g)$ above admits a homogeneous Riemannian structure if and
only if it is a reductive homogeneous Riemannian manifold.  This means that
$M=G/H$, where $G$ is a connected Lie group acting transitively and
effectively on~$M$ via isometries, $H$~is the isotropy group at a base
point $o\in M$, and the Lie algebra~$\lie g$ of~$G$ may be decomposed into
a vector space direct sum $\lie g=\lie h+\lie m$, where $\lie h$~is
the Lie algebra of~$H$ and $\lie m$ is an $\Ad(H)$-invariant subspace,
i.e., $\Ad(H)\lie m\subset\lie m$.  As $G$ is connected and $M$
simply-connected, $H$ is connected, and the latter condition is equivalent to
$[\lie h,\lie m] \subset \lie m$.

Conversely, let $S$ be a homogeneous Riemannian structure on a complete
Riemannian manifold~$(M,g)$ that is connected and simply-connected.  We fix
a point $o\in M$ and put $\lie m=T_oM$.  If $\Rt$ is the curvature tensor
of $\Nt$, we can consider the holonomy algebra $\tilde{\lie h}$ of $\Nt$ as
the Lie subalgebra of skew-symmetric endomorphisms of $(\lie m,g_o)$
generated by the operators $\Rt_{XY}$, where $X,Y\in\lie m$.  Then,
according to the Ambrose-Singer construction
\cite{Ambrose-S:homogeneous,Tricerri-V:homogeneous}, a Lie bracket is
defined in the vector space direct sum $\tilde{\lie g}=\tilde{\lie
h}+\lie m$ by
\begin{equation}
  \label{corchete}
  \left\{
    \begin{alignedat}{2}
      [U,V] &= UV-VU,  &&U,V\in\tilde{\lie h}, \\
      [U,X] &= U(X),   &&U\in\tilde{\lie h},\ X\in\lie m, \\
      [X,Y] &= S_XY - S_YX + \Rt_{XY}, &\qquad&X,Y\in\lie m.
  \end{alignedat}
  \right.
\end{equation}
One calls $(\tilde{\lie g},\tilde{\lie h})$ the reductive pair associated
to the homogeneous Riemannian structure~$S$.  The connected,
simply-connected Lie group $\widetilde{G}$ whose Lie algebra
is~$\tilde{\lie g}$ acts transitively on~$M$ via isometries and $M\equiv
\widetilde{G}/\widetilde{H}$, where $\widetilde{H}$ is the connected Lie
subgroup of~$\widetilde{G}$ whose Lie algebra is~$\tilde{\lie h}$.  The set
$\Gamma$ of elements of $\widetilde{G}$ which act trivially on~$M$ is a
discrete normal subgroup of $\widetilde{G}$, and the Lie group
$G=\widetilde{G}/\Gamma$ acts transitively and effectively on $M$ as a
group of isometries, with isotropy group $H=\widetilde{H}/\Gamma$. Then,
there exists a diffeomorphism $\varphi \colon G/H \to M$ and $(M,g)$ is
(isometric to) the reductive homogeneous Riemannian manifold $(G/H,\varphi^*g)$.

\subsection{Homogeneous quaternionic K\"ahler structures}

We recall that an almost quaternionic structure on a $C^\infty$
manifold~$M$ is a rank~$3$ subbundle $\upsilon$ of the bundle of $(1,1)$
tensors on~$M$, such that there locally exists a basis $J_1$, $J_2$, $J_3$
satisfying the conditions
\begin{equation}
  \label{jujdjt}
  J^2_1 = J^2_2 = J^2_3 = -I, \quad
  J_1J_2 = -J_2J_1 = J_3, \quad \text{etc.}
\end{equation}
Here and throughout the rest of this paper, `etc.' denotes the equations
obtained by cyclically permuting the indices.

Such a basis is called a standard local basis of~$\upsilon$ in its domain
of definition.  Then, $(M, \upsilon)$ is called an almost quaternionic
manifold, and $M$ has dimension $4n$, with $n \geqslant1$.  On any almost
quaternionic manifold $(M, \upsilon)$, there is a Riemannian metric $g$
such that $g(\sigma X, Y) + g(X, \sigma Y) = 0$, for any section $\sigma$
of $\upsilon$.  Then, $(M,g,\upsilon)$ is called an almost
quaternion-Hermitian manifold.  It is known that $M$ admits an almost
quaternion-Hermitian structure if and only if the structure group of the
tangent bundle~$TM$ is reducible to~$\SP(n)\SP(1)$
(cf.~\S\ref{sec:action}).

Let $J_1, J_2, J_3$ be a standard local basis of~$\upsilon$ and let
\begin{equation*}
  \omega_a (X, Y) = g(X, J_a Y),\qquad a=1,2,3.
\end{equation*}
These are local differential forms, but the differential $4$-form $\Omega =
\sum_{a=1}^3 \omega_a \wedge\omega_a$ is known to be globally defined.  Note that we
have
\begin{equation}
  \label{jhjh}
  g( J_aX, J_aY ) = g(X,Y), \qquad a=1,2,3.
\end{equation}

The manifold is said to be quaternionic K\"ahler if $\LC\Omega=0$ or,
equivalently, one has locally \cite{Ishihara:qK} that
\begin{equation}
  \label{naji}
  \LC_X J_1 = \tau^3(X)J_2 - \tau^2(X)J_3, \quad \text{etc.},
\end{equation}
for certain differential $1$-forms $\tau^1,\tau^2,\tau^3$.

In the present paper we shall consider quaternionic K\"ahler manifolds of
$\dim\geqslant 8$ and non-zero scalar curvature (see
\cite{Besse:Einstein,Salamon:holonomy,Swann:MathAnn}).

\begin{definition} (\cite[p.\ 218]{Alekseevsky-C:isometry}) A quaternionic
  K\"ahler manifold $(M,g,\upsilon)$ is said to be a \emph{homogeneous
  quaternionic K\"ahler manifold\/} if it admits a transitive group of
  isometries.
\end{definition}

\begin{remark}
  Concerning isometry groups, the situation for quaternionic K\"ahler
  manifolds is rather different from that for K\"ahler (see
  \cite[p.~375]{Abbena-G:aH}).  In fact, a quaternionic K\"ahler
  manifold~$M$ with $\dim M \geqslant8$ and non-zero scalar curvature is,
  even locally, irreducible \cite{Berger:hol}, and its Ricci tensor is
  nowhere zero.  Thus, by a theorem of Kostant \cite{Kostant:holonomy}, a
  transitive group of isometries induces the Lie algebra of the restricted
  holonomy group, and thus preserves $\Span\{ J_1,J_2,J_3 \}$, since each
  $J_a$ belongs to the Lie algebra of the holonomy group (see
  \cite[p.~407]{Besse:Einstein}).
\end{remark}

On the other hand, we have the following Corollary of Kiri\v cenco's
Theorem \cite{Kiricenko:homogeneous} (see also
\cite{Ambrose-S:homogeneous,Fino:torsion}).

\begin{theorem}
  \label{thqKhm}
  A connected, simply-connected and complete quaternionic K\"ahler manifold
  $(M, g, \upsilon)$ is homogeneous if and only if there exists a tensor
  field $S$ of type $(1,2)$ on $M$ satisfying
  \begin{equation*}
    \Nt g=0, \quad\Nt R=0, \quad\Nt S=0, \quad\Nt\Omega= 0,
  \end{equation*}
  where $\Nt=\LC-S$.
\end{theorem}

Such a tensor~$S$ is called a \emph{homogeneous quaternionic K\"ahler
structure} on~$M$.

The equation $\Nt\Omega=0$ is equivalent, under $\Nt g = 0$, to the
existence of three differential $1$-forms $\tilde\tau^1$, $\tilde\tau^2$,
$\tilde\tau^3$ such that
\begin{equation}
   \label{nabj}
   \Nt_X J_1 =\tilde\tau^3(X)J_2 - \tilde\tau^2(X)J_3, \quad \text{etc.}
\end{equation}
Combined with~\eqref{naji}, the previous formul\ae\ yield
\begin{equation}
  \label{spospnt}
  S_X J_1Y-J_1S_X Y = \pi^3(X)J_2Y-\pi^2(X)J_3Y,\quad \text{etc.},
\end{equation}
for $\pi^a =\tau^a -\tilde\tau^a$, $a=1,2,3$.  Writing as usual $S_{XYZ} =
g(S_XY, Z)$, we have that
\begin{equation}
  \label{sjiji}
  S_{XJ_1YJ_1Z} - S_{XYZ}  = \pi^3(X)g( J_2Y, J_1Z) -\pi ^2(X)g(
  J_3Y,J_1Z), \quad \text{etc.,}
\end{equation}
which together with the condition $S_{XYZ} = -S_{XZY}$, are the symmetries
satisfied by a homogeneous quaternionic K\"ahler structure~$S$.

Note moreover that~$S_X$ acts as an element of the Lie algebra $\sP(1)
\oplus \sP(n)$ on~$T_pM$, for any $p \in M$.  In fact, from the definition
of~$\SP(n)\SP(1)$ (see \eqref{spnspot} below), an element $U \in \sP(1)
\oplus \sP(n)$ is characterised by the condition $U \circ J_a - J_a \circ U
= m^b_a J_b$, for a matrix $(m^b_a) \in \so(3)$.
Equation~\eqref{spospnt} shows that this is indeed satisfied by~$S$.

\subsection{Fino's classification}
\label{finor}

Let $E$ denote the standard representation of $\SP(n)$ on $\mathbb C^{2n}$.
This representation is quaternionic, meaning that it carries an anti-linear
endomorphism~$j$ that commutes with the action of~$\SP(n)$ and
satisfies~$j^2=-1$.  Write $S^rE$ for the $r$th-symmetric power of~$E$, so
$S^2E\cong\sP(n)\otimes\mathbb C$, and let $K$~be the irreducible
$\SP(n)$-module in $E\otimes S^2E = S^3E + K + E$, ($K$~is of highest weight
$(2,1,0, \dots, 0)$).  Take $H$ to be the standard representation of
$\SP(1)\cong\SU(2)$ on $\mathbb C ^2$, then $S^2H\cong\sP(1)\otimes\mathbb
C$ and $S^3H$ is the $4$-dimensional irreducible representation
of~$\SP(1)$.

Homogeneous quaternionic K\"ahler structures are classified from a
representation-theoretic point of view as follows.

\begin{theorem}[Fino~{\cite[Lemma 5.1]{Fino:torsion}}]
  \begin{equation*}
    \begin{split}
      \mathcal T(V)_+ &= [EH] \otimes(\sP(1) \oplus \sP(n))\\
      &\cong [EH] + [ES^3H] + [EH] + [S^3EH] + [KH].
    \end{split}
  \end{equation*}
\end{theorem}

Here, $[V]$ denotes the real representation whose complexification is~$V$,
sums are direct, and the tensor products signs are omitted, that is, one
writes $EH$ instead of $E \otimes H$, and so on.  We shall write $\QK1$,
\dots, $\QK5$ for the five Fino classes in the above order, which differs
slightly from Fino's.  Thus $\QK1=[EH]\subset[EH]\otimes\sP(1)$, etc.  We
also write $\QK{i+j}$ for $\QK i+\QK j$, etc.

\subsection{Some conventions}

We shall use the following conventions for the curvature tensor of a linear
connection of a Riemannian manifold $(M,g)$:
\begin{gather*}
  R_{XY}Z=\LC_{[X,Y]}Z-\LC_X\LC_Y Z + \LC_Y\LC_X Z, \\
  R_{XYZW} = g(R_{XY}Z,W), \quad R_{XY}(Z,W) = R_{XYZW}.
\end{gather*}
We denote the Ricci tensor by $\mathbf{r}$ and the scalar curvature by
$\mathbf{s}$.  We write $\nu = \mathbf{s}/4n(n+2)$ for the reduced scalar
curvature of a $4n$-dimensional Riemannian manifold.  In addition, the
Einstein summation convention for repeated indices is assumed.

\section{Classification by real tensors}
\label{sec:cla}

\subsection{The space of tensors}

Let $(V, \inp\cdot\cdot, J_1, J_2, J_3)$ be a quaternion-Hermitian real
vector space, i.e., a $4n$-dimensional real vector space endowed with an
inner product $\inp\cdot\cdot$ and operators $J_1,J_2,J_3$ satisfying
\eqref{jujdjt} and \eqref{jhjh}.  Such a space $V$ is the model for the
tangent space at any point of a quaternionic K\"ahler manifold.  Consider
the space of tensors
\begin{equation*}
  \mathcal{T}(V)=\{\, S\in\otimes^3V^*:S_{XYZ}=-S_{XZY} \,\}
\end{equation*}
and the vector subspace $\V$ of $\mathcal{T}(V)$ defined by
\begin{equation*}
  \V = \{\, S\in \otimes^3V^*:S_{XYZ}=-S_{XZY},\
  \exists \pi^a \in V^*\ \text{s.t. $S$~satisfies~\eqref{sjiji}} \,\}.
\end{equation*}
Any homogeneous Riemannian structure on~$M$ belongs to $\mathcal T(T_pM)$
pointwise, whereas homogeneous quaternionic K\"ahler structures are
pointwise in~$\V$.  We wish to explicitly decompose~$\V$.

For each element $S\in\V$ consider the tensor
\begin{equation}
  \Theta_{XYZ}^S = \tfrac12\pi^a(X)\inp{J_a Y}Z, \label{fT}
\end{equation}
which up to a factor $-4$ is the sum of the right-hand sides
of~\eqref{sjiji}.

\begin{lemma}
  Given $S\in\V$, the tensor $\Theta^S$ lies in~$\V$ and satisfies the
  equalities \eqref{sjiji} with the same forms $\pi^a$, $a=1,2,3$, as~$S$.
\end{lemma}

\begin{proof}
  This follows directly from the relations $J_1 J_2 = J_3$, etc., and
  \eqref{jhjh}.
\end{proof}

On account of the equalities \eqref{sjiji}, for any $S\in\V$ we have that
\begin{equation*}
  S_{XYZ}=\Theta_{XYZ}^S+T^S_{XYZ},
\end{equation*}
where
\begin{equation}
  T^S_{XYZ}=\frac14\bigl(S_{XYZ}+ \sum_{a=1}^3S_{XJ_aYJ_aZ}\bigr).
\label{fTT}
\end{equation}
The tensor $T^S$ belongs to
\begin{equation*}
  \VE=\{\,T\in\otimes^3V^*:T_{XYZ}=-T_{XZY}, \ T_{XJ_a YJ_a Z}=T_{XYZ}\
  \forall a\,\},
\end{equation*}
that is, $\VE$ is the subspace of $\V$ defined by the conditions $\pi^a
=0$.  We also define, corresponding to $T = 0$, the subspace of $\V$
\begin{equation*}
  \VH=\{\,\Theta\in\otimes^3V^*:\Theta_{XYZ}=\tfrac12
  \pi^a (X)\inp{J_a Y}Z,\ \pi^a  \in V^*\, \},
\end{equation*}
which can be also given as
\begin{equation}
  \label{otroch}
  \VH=\{\,S\in\V: S_{XYZ}+ \sum_{a=1}^3 S_{XJ_aYJ_aZ}= 0\,\}.
\end{equation}

\begin{proposition}
  The space $\V$ decomposes as an orthogonal direct sum
  \begin{equation}
    \V=\VH + \VE \label{fD1}
  \end{equation}
  with respect to the inner product
  \begin{equation}
    \inp  S{S'}  = \sum_{r,s,t=1}^{4n} S_{e_re_se_t} S'_{e_re_se_t}, \label{fPE}
  \end{equation}
  where $\{e_r\}_{r=1, \dots, 4n}$ is any orthonormal basis of $V$.
\end{proposition}

\begin{proof}
  If $S \in \V$, we have already seen that we can write it as $S=\Theta+T$,
  with $\Theta$ and $T$ defined in \eqref{fT} and~\eqref{fTT}.  Conversely,
  put $S = \Theta+T$ with $\Theta_{XYZ} = \frac12 \pi^a (X) \inp{J_a Y}Z$
  for some one-forms $\pi^a$ and $T\in\VE$.  It is easily checked that $S$
  satisfies \eqref{sjiji} for the forms $\pi^a$, so $S\in\V$.  To prove
  that the decomposition \eqref{fD1} is orthogonal, we take an orthonormal
  basis of~$V$ of the form $\{e_r\}=\{u_s,J_1u_s,J_2u_s,J_3
  u_s\}_{s=1,\dots,n}$.  Then, for $T\in\VE$ and $\Theta \in\VH$, we have
  that
  \begin{equation*}
    \inp T\Theta =-\frac14\sum_{r=1}^{4n}\sum_{s=1}^n \pi^a(e_r)
    \Bigl(T_{e_r J_a u_su_s}-\sum _{b=1}^3T_{e_r J_b J_a J_b u_su_s}\Bigr)=0.
\end{equation*}
\end{proof}

\subsection{The action of $\SP(n)\SP(1)$}
\label{sec:action}

Standard local bases of~$\upsilon$ on a quaternionic K\"ahler manifold
$(M,\upsilon)$ are not intrinsic.  In fact, given one basis
$\{J_1,J_2,J_3\}$, the other bases $\{J_1',J_2',J_3'\}$ are obtained as
\begin{equation*}
  J_a' = m_a^b J_b , \qquad a,b = 1,2,3,
\end{equation*}
for arbitrary $(m_b^a)\in \SO(3)$.  On the other hand, by its very
definition, the holonomy group of a quaternionic K\"ahler manifold is
contained in
\begin{equation*}
  \SP(n)\SP(1) = (\SP(n)\times \SP(1))/\{\pm \Id\}\subset \SO(4n).
\end{equation*}

The action of this group on $\mathbb R^{4n} \equiv \mathbb H^n$ is as follows:
\begin{equation}
  \label{srspnspud}
  (B,q)v = B v \bar q, \qquad B \in \SP(n), \ q \in \SP(1),
\end{equation}
where the $v$ are regarded as vectors in $\mathbb H^n$ and $\bar q$ denotes
the quaternionic conjugate of $q$.  It is easy to check that an orthogonal
automorphism $A \in \SO(4n)$ belongs to $\SP(n)\SP(1)$ if and only if
\begin{equation}
  \label{spnspot}
  A \circ J_a = m^b_a J_b \circ A,
\end{equation}
for a certain matrix $(m_b ^a )\in \SO(3)$, which is obtained from the
projection homomorphism
\begin{equation*}
  \SP(n)\SP(1) \longrightarrow \SP(1)/\{\pm \Id\}=\SO(3).
\end{equation*}

The standard representation of $\SP(n)\SP(1)$ on $V$, defined by
\eqref{srspnspud}, induces a representation of $\SP(n)\SP(1)$ on $\V$ given
by
\begin{equation}
  \label{asxyz}
  (A(S))_{XYZ} = S_{A^{-1}XA^{-1}YA^{-1}Z}.
\end{equation}

\begin{proposition}
  The subspaces $\VH\subset\V$ and $\VE\subset\V$ are invariant under the
  action of $\SP(n)\SP(1)$ on $\V$.
\end{proposition}

\begin{proof}
  Let $\Theta=\frac12\pi^a (X)\inp{J_a Y}Z \in\VH$.  Then
  \begin{equation*}
    (A(\Theta))_{XYZ} = \Theta_{A^{-1}X A^{-1}YA^{-1}Z}
    = \frac12\pi^a(A^{-1}X)\inp{J_a A^{-1}Y}{A^{-1}Z}.
  \end{equation*}
  Since $A J_a = m_a^b J_b A$, for $a=1,2,3$, we have that
  \begin{equation*}
    (A(\Theta))_{XYZ}
    = \tfrac12 \pi^a (A^{-1}X)\inp{m_a^b A^{-1} J_b Y}{A^{-1}Z}
    = \tfrac12 \bar\pi^b (A^{-1}X)\inp{J_b Y}Z,
  \end{equation*}
  because $A\in \SO(4n)$ and we put $\bar\pi^b = m_a^b \pi^a$.  So
  $A(\Theta)$ belongs to $\VH$, with the forms $\bar{\pi}^b\circ A^{-1}$.

  Let $T\in\VE$.  Then we have that
  \begin{equation}
    \label{fII}
    T_{XJ_a YJ_a Z}=T_{XYZ},\qquad a=1,2,3.
  \end{equation}
  Then $(A(T))_{XJ_a YJ_a Z}=(A(T))_{XYZ}$.  Indeed, from \eqref{fII} we
  have that
  \begin{equation*}
    \begin{split}
      (A(T))_{XJ_a YJ_a Z}
      &= m_a^b m_a^c T_{A^{-1}X J_b A^{-1}Y J_c A^{-1}Z}\\
      &= - m_a^b m_a^c T_{A^{-1}X A^{-1}Y J_bJ_c A^{-1}Z} \\
      &= \Bigl(\sum_{b=1}^3 (m_a^b)^2\Bigr)
      T_{A^{-1}XA^{-1}YA^{-1}Z}=(A(T))_{XYZ},
    \end{split}
  \end{equation*}
  when $(m_a^b)$ is the matrix associated to~$A^{-1}$.
\end{proof}

We conclude that the representation~\eqref{asxyz} of $\SP(n)\SP(1)$
decomposes as $\V = \VH + \VE$.  However, neither space is
irreducible.  We devote the next two subsections to the explicit
decomposition of each space in turn.

\subsection{The space $\VH$}

In the space $\VH$ one can first distinguish the subspace
\begin{equation*}
  \VH_0^\perp=\{\, \Theta\in\VH :
  \Theta _{XYZ}= \sum_{a=1}^3 \theta(J_aX)\inp{J_aY}Z ,\ \theta\in V^* \,\}.
\end{equation*}

\begin{proposition}
  The space $\VH_0^\perp$ is $\SP(n)\SP(1)$-invariant.
\end{proposition}

\begin{proof}
  Given $A\in \SP(n)\SP(1)$, we have that
  \begin{equation*}
    \begin{split}
      (A(\Theta) )_{XYZ}
      &=\sum _{a=1}^3\theta(J_a (A^{-1}X))\inp{J_a A^{-1}Y}{A^{-1}Z} \\
      &=\sum_{a=1}^3 m_a^b m_a ^c \theta(A^{-1}J_b X)\inp{A^{-1} J_c
      Y}{A^{-1}Z}\\
      &= \sum _{a=1}^3 \theta(A^{-1}J_a X)\inp{J_a Y}Z.
    \end{split}
  \end{equation*}
  So $A(\Theta) $ belongs to $\VH_0^\perp$ with the form $\theta\circ
  A^{-1}$.
\end{proof}

We now describe the orthogonal subspace in $\VH$.

\begin{proposition}
  The space orthogonal to $\VH_0^\perp$ in $\VH$ with respect to the scalar
  product defined in~\eqref{fPE}, is
  \begin{equation*}
    \VH_0=\{\,\Theta\in\VH:\Theta_{XYZ}=\pi ^a(X)\inp{J_a Y}Z ,\ \pi^a \circ
    J_a =0\,\}.
  \end{equation*}
  This subspace is invariant under the action of $\SP(n)\SP(1)$.  So $\VH=
  \VH_0^\perp + \VH_0$ as $\SP(n)\SP(1)$-modules.
\end{proposition}

\begin{proof}
  Let $\{e_r\}$ be an orthonormal basis of~$V$ as above.  Then, for
  $\bar{\Theta}\in\VH_0^\perp$ and $\Theta\in \VH$, the condition
  $\inp*{\bar{\Theta}}\Theta =0$ gives us
  \begin{equation*}
    \begin{split}
      0 &=\sum_{r,s,t=1}^{4n}\bar{\Theta}_{e_re_se_t}\Theta_{e_r e_se_t} =
      4n
      \sum_{r=1}^{4n} \theta(J_a e_r)\pi^a (e_r)\\
      & =-4n\sum_{r=1}^{4n} \theta(e_r)\pi^a (J_a e_r),
    \end{split}
  \end{equation*}
  and this happens for any form $\theta$ if and only if $\pi^a \circ J_a
  =0$.
\end{proof}

We now give another characterisation of $\VH_0$.

\begin{proposition}
  For $\dim V=4n$ with $n>1$, one has
  \begin{equation*}
    \VH_0= \biggl\{\,\Theta \in \VH:
    \Cyclic_{XYZ} \Theta_{XYZ} + \sum_{a=1}^3
    \Cyclic_{XJ_aYJ_aZ} \Theta_{XJ_aYJ_aZ}=0 \,\biggr\}.
  \end{equation*}
\end{proposition}

\begin{proof}
  It is not difficult to check that if $\Theta_{XYZ} = \pi^a (X)\inp{J_a
  Y}{Z } \in\VH$, then one has
  \begin{equation*}
    \begin{split}
      {\Cyclic_{XYZ}} &\Theta_{XYZ} +\sum_{a=1}^3{\Cyclic_{XJ_aYJ_aZ}}
      \Theta_{XJ_aYJ_aZ} \\
      &=(\pi^a \circ J_a )\bigl( \inp XY Z -\inp XZ Y \bigr) \\
      &\qquad +{\Cyclic_{123}} (\pi^1\circ J_2 - \pi^2\circ
      J_1+\pi^3)\bigl(\inp{J_3 Z}X Y -\inp{ J_3Y}X Z\bigr).
    \end{split}
  \end{equation*}

  If $\Theta \in \VH_0$, then the right-hand side of the previous equation
  vanishes.  Conversely, for $X=Z$ orthogonal to $\Span \{Y, J_1Y, J_2Y,
  J_3Y \}$, one obtains that $(\pi^a \circ J_a )(Y)=0$ for all $Y\in V$.
\end{proof}

Let
\begin{equation*}
  c_{12}(S)(Z) = \sum_{r=1}^{4n} S_{e_re_rZ}
\end{equation*}
for any orthonormal basis $\{ e_r \}$, $r =1,\dots,4n$.  Then we have the
next characterisation of $\VH_0$, which justifies the notation.

\begin{proposition}
  We have that $\VH_0 = \{\Theta\in\VH : c_{12}(\Theta)=0\}$.
\end{proposition}

\begin{proof}
  Let $\{e_r\}$ be an orthonormal basis of $V$ as above.  For $\Theta \in
  \VH$ and $t\in\{1,\dots,4n\}$ fixed, we have that
  \begin{equation*}
    c_{12}(\Theta)(e_t)=\sum_{s=1}^{4n}\pi^a
    (e_s)\inp{J_a e_s}{e_t} =-\pi^a (J_a
    e_t)\inp{e_t }{e_t} =-(\pi^a \circ J_a
    )(e_t).
\end{equation*}
\end{proof}

An alternative characterisation of~$\VH_0$, based on the
expression~\eqref{otroch} of the elements of $\VH$, is given by the next
proposition.

\begin{proposition}
  We have that
  \begin{equation*}
    \VH_0=\Bigl\{\,\Theta\in\VH:\Theta_{XYZ}+
    \sum_{a=1}^3\Theta_{J_aXJ_aYZ}=\Theta_{YXZ}+
    \sum_{a=1}^3\Theta _{J_aYJ_aXZ}\,\Bigr\}.
  \end{equation*}
\end{proposition}

\subsection{The space $\VE$}

Consider the map $L \colon \VE \rightarrow \VE$ defined by
\begin{equation*}
  L(T)_{XYZ}= T_{ZXY}+T_{YZX}+ \sum_{a=1}^3\bigl(T_{J_aZXJ_aY}+T_{J_aYJ_aZX}\bigr).
\end{equation*}
It is easily seen that $L(T)\in \VE$, for any $T$ and that $L$~is a linear
map satisfying $L(A(T))=A(L(T))$.  Moreover, we have the next results.

\begin{proposition}
  The map $L$ satisfies $L\circ L = 8\Id-2L$. \qed
\end{proposition}

\begin{corollary}
  The minimal polynomial of $L$ is $(x-2)(x+4)$.  Thus, $L$ is
  diagonalisable with two eigenspaces $\VE^2$ and $\VE^{-4}$ with
  respective eigenvalues $2$ and $-4$, and $\VE=\VE ^2 + \VE^{-4}$.
  Since $L(A(T))=A(L(T))$, these eigenspaces are invariant under the action
  of $\SP(n)\SP(1)$. \qed
\end{corollary}

\begin{proposition}
  The subspaces $\VE^2$ and $\VE^{-4}$ are mutually orthogonal.
\end{proposition}

\begin{proof}
  A straightforward calculation shows that $L$ is self-adjoint; that is,
  that $\inp{L(T)}{T'} = \inp T{L(T')}$, for $T,T' \in \VE$.  Then, taking
  $T \in\VE^2$ and $T' \in\VE^{-4}$, one obtains that $2\inp T{T'} = -4\inp
  T{T'}$, thus concluding.
\end{proof}

\subsubsection{The subspace $\VE^2$}

A tensor $T$ belongs to $\VE^2$ if and only if $L(T) = 2T$, so
\begin{equation}
  \label{unsexto}
  \VE^2 = \biggl\{\,T\in\VE : T_{XYZ}=\frac 16\Bigl({\Cyclic_{XYZ}} T_{XYZ}
  + \sum_{a=1}^3{\Cyclic_{XJ_aYJ_aZ}} T_{XJ_aYJ_aZ}\Bigr)\,\biggl\}.
\end{equation}
For $\theta\in V^*$, put
\begin{equation*}
  \begin{split}
    T^\theta_{XYZ}
    &= \inp XY \theta(Z)-\inp XZ\theta(Y) \\
    &\qquad + \sum_{a=1}^3 \bigl(\inp X{J_aY} \theta(J_aZ) -\inp
    X{J_aZ}\theta(J_aY)\bigr).
  \end{split}
\end{equation*}
From the expression of the tensors in the subspace $\VH_0^\perp$, we now
consider the space
\begin{equation*}
    \VE_0^\perp = \bigl\{\, T^\theta \in \VE : \theta\in V^*\,\bigr\}.
\end{equation*}

The tensors $T^\theta$ in $\VE_0^\perp$ satisfy the cyclic sum property
\eqref{unsexto}, showing $\VE_0 ^\perp \subset \VE ^2$, and also the
condition $c_{12}(T^\theta)=4(n+1)\theta$.  It is also straightforward to
check the invariance of $\VE_0^\perp$ under the action of $\SP(n)\SP(1)$.
Indeed, $A(T^\theta)=T^{\theta\circ A^{-1}}$, for $A \in \SP(n)\SP(1)$.

\begin{proposition}
  The subspace orthogonal to $\VE_0^\perp$ in $\VE$ is the subspace defined
  by $\VE_0=\{T\in\VE:c_{12}(T)=0\}$.  This space is thus invariant under
  the action of $\SP(n)\SP(1)$.
\end{proposition}

\begin{proof}
  Let $\{e_r\}$ be an orthonormal basis of $V$ as above.  Let $T \in\VE_0$.
  Then $\inp{T^\theta}T=0$ for any $T^\theta \in\VE_0^\perp$.  In particular,
  taking $\theta$ to be the dual basis element to~$e_\ell$,
  $\ell\in\{1,\dots,4n\}$, we have that
  \begin{equation*}
    \inp{T^\theta}T = 2\sum_{r=1}^{4n}(T_{e_re_re_l}+ T_{J e_rJ e_re_l}) =
    8c_{12}(T)(e_l).
  \end{equation*}
\end{proof}

\subsubsection{The subspace $\VE^{-4}$}

A tensor $T$ belongs to $\VE^{-4}$ if and only if $L(T) = -4T$, so
\begin{equation}
  \label{fV2}
  \VE^{-4} =\bigl\{\, T \in \VE :{\Cyclic_{XYZ}} T_{XYZ}+
  \sum_{a=1}^3{\Cyclic_{XJ_aYJ_aZ}} T_{XJ_aYJ_aZ}=0\,\bigr\} .
\end{equation}

\begin{proposition}
  One has $\VE^{-4} = \bigl\{T\in\VE : \Cyclic_{XYZ} T_{XYZ}=0\bigr\}$.
\end{proposition}

\begin{proof}
  This is immediate from
  \begin{equation*}
    0 = \Cyclic_{XYZ} \bigl( \Cyclic_{XYZ} T_{XYZ} +
    \sum _{a=1}^3\Cyclic_{XJ_aYJ_aZ} T_{XJ_aYJ_aZ}\bigr) = 4\Cyclic_{XYZ}T_{XYZ}.
  \end{equation*}
\end{proof}

\begin{proposition}
  The space $\VE^{-4}$ is contained in $\VE_0$, so the tensors
  $T\in\VE^{-4}$ are traceless with respect to~$c_{12}$.
\end{proposition}

\begin{proof}
  Let $\{e_r\}$ be an orthonormal basis of $V$ as above.  If $T\in
  \VE^{-4}$, then it satisfies the cyclic sum condition as in~\eqref{fV2},
  and in particular we have that
  \begin{equation*}
    0 = 3\,T_{e_se_se_t}+ \sum _{a=1}^3(T_{J_ae_sJ _a e_se_t }+T_{J_a
    e_te_sJ_a e_s}).
  \end{equation*}
  Summing over~$s$, one gets $0=6c_{12}(T)(e_t)+ \sum _{a=1}^3
  \sum_{s=1}^{4n} T_{J_a e_te_sJ_a e_s}$.

  We claim that the second summand vanishes.  Indeed, by the cyclic sum
  property we have that
  \begin{equation*}
    \begin{split}
      0 & =\sum_{a=1}^3 \sum_{s=1}^{4n} \Bigl(T_{J_a e_te_sJ_ae_s} +
      T_{e_sJ_a e_sJ_a e_t} +T_{J_a e_sJ_a e_te_s} \\
      & \qquad\qquad + \sum_{b=1}^3(T_{J_a e_tJ_b e_sJ_b J_a e_s} +
      T_{J_be_sJ_b J_a
      e_sJ_a e_t} + T_{J_b J_a e_sJ_a e_t J_b e_s} )\Bigr) \\
      & =4\sum_{a=1}^3\sum_{s=1}^{4n} \Bigl(T_{J_a e_te_sJ_ae_s} +
      T_{e_se_se_t} + T_{J_ae_sJ_a e_se_t}\\
      &\qquad\qquad\qquad -\sum _{b=1}^3( T_{J_b e_sJ_a J_b J_a e_se_t } +
      T_{J_bJ_ae_sJ_bJ_ae_te_s})\Bigr).
    \end{split}
  \end{equation*}
  Evaluating the last four terms we obtain multiples of~$c_{12}(T)(e_t)$
  that cancel and so we are left with the first term being zero, from which
  the result follows.
\end{proof}

As $\VE^2$ is orthogonal to $\VE^{-4}$, then $\VE^2\cap \VE_0$ is another
subrepresentation, orthogonal to the above ones.  According to the previous
results, one has the orthogonal decomposition
\begin{equation*}
  \V =
  \VH_0^\perp  + \VH_0
  +  \VE_0^\perp   +
  (\VE^2  \cap \VE_0)  +
  \VE^{-4}.
\end{equation*}

\subsection{The classification theorem}

\begin{theorem}
  \label{tchqKs}
  If $n \geqslant2$, then $\V$ decomposes into the direct sum of the
  following subspaces invariant and irreducible under the action of
  $\SP(n)\SP(1)$:
  \begin{align*}
    \QK1 & = \{\Theta\in\VH:\Theta_{XYZ} = \sum_{a=1}^3\theta(J_a X)
    \inp{J_a Y}Z ,\ \theta\in V^*\},\\
    \QK2 & =\bigl\{\Theta\in\VH:\Theta_{XYZ} = \theta^a (X) \inp{J_a Y}Z,\
    \theta^a \circ J_a =0,\\
    &\hspace{18em}\theta^1, \theta^2, \theta^3\in V^* \bigr\}, \\
    \QK3 & = \{T \in \VE : T_{XYZ}
    = \inp XY\theta(Z)-\inp XZ\theta(Y)\\
    &\hspace{5em} + \sum_{a=1}^3\bigl( \inp X{J_a Y} \theta(J_a Z)-
    \inp X{J_a Z }
    \theta(J_a Y)\bigr),\theta\in V^*\}, \\
    \QK4 & =\bigl\{ T \in \VE \colon T_{XYZ} = \frac16 \bigl( \Cyclic_{XYZ}
    T_{XYZ} + \sum\limits_{a=1}^3
    \Cyclic_{XJ_aYJ_aZ} T_{XJ_aYJ_aZ} \bigr),\\
    &\hspace{20em} \sum_{r=1}^{4n}T_{e_re_rZ}=0 \bigr\},\\
    \QK5 & = \bigl\{ T \in \VE : \Cyclic_{XYZ} T_{XYZ} = 0 \bigr\}.
  \end{align*}
\end{theorem}

In other words, $\QK1=\VH_0^\perp$, $\QK2=\VH_0$, $\QK3=\VE_0^\perp$,
$\QK4=\VE^2\cap\VE_0$ and $\QK5=\VE^{-4}$.

\begin{proof}
  Noting that $\V\cong[EH]\cong\V^*$, it suffices to identify the five
  modules above with the modules in Fino's classification (see
  \S\ref{finor})
  \begin{equation*}
    [EH] +  [ES^3H]  +  [EH] + [S^3EH] +  [KH].
  \end{equation*}

  It is clear that $\QK1 \cong [EH] \cong \QK3$, with $\QK1$ being the copy
  in $[EH] \otimes [S^2H]$ and $\QK3$ the copy in $[EH] \otimes [S^2E]$.
  The definition of $\QK2$ shows that it is the complement of $[EH] \cong
  \QK1$ in $[EH] \otimes [S^2H]$, so $\QK 2 \cong [ES^3H]$.

  Now we have that
  \begin{equation*}
    \QK{3+4+5} \cong [EH] \otimes [S^2E] \cong [EH] + [S^3EH] + [KH],
  \end{equation*}
  and know that $\QK3 \cong [EH]$.

  We can identify the module $[KH]$ as follows.  Note that $[KH]$ is both a
  submodule of $[EH] \otimes [S^2E] \subset V^*\otimes \Lambda^2V^*$ and of
  $[\Lambda^2E] \otimes [EH] \subset S^2V^* \otimes V^*$.  Indeed,
  \begin{equation*}
    [\Lambda^2E] \otimes [EH] \cong 2[EH] + [\Lambda^3_0EH] + [KH],
  \end{equation*}
  which does not contain $[S^3EH]$.  Thus an equivariant map $\QK{4+5} \to
  [\Lambda^2E] \otimes [EH]$ will always contain $[S^3EH]$ in its kernel
  and can only be non-zero on a module isomorphic to $[KH]$.  The module
  $[\Lambda^2E] \subset S^2V^*$ consists of symmetric bilinear forms $b$
  with $b(J\cdot,J\cdot)= b$ for each $J$, so let $p\colon S^2V^* \otimes
  V^* \to [\Lambda^2E] \otimes [EH]$ be given by
  \begin{equation*}
    p(T)_{XYZ} = \frac14 \Bigl(T_{XYZ}+ \sum _{a=1}^3T_{J_aX J_aY Z}\Bigr).
  \end{equation*}

  Consider $\QK5$, which consists of the $T \in V^*\otimes \Lambda^2V^*$
  such that $T_{XYZ}+T_{YZX}+T_{ZXY}=0$.  The projection to this module is
  given by mapping $T$ to $U$, where $U_{XYZ} = \frac14(T_{XYZ} + \sum
  _{a=1}^3T_{XJ_aYJ_aZ})$ and then by mapping $U$ to~$\frac16(2-L)U$.

  Applying these maps to the element $\alpha\otimes \beta\wedge\gamma \in
  V^*\otimes \Lambda^2V^*$, the projection to $\QK5$ is
  \begin{multline*}
    \tfrac1{12}\bigl(\ow\A\B\Cc-2\ow\B\Cc\A-2\ow\Cc\A\B\\
    + \sum _{a=1}^3(\ow\A{J_a\B}{J_a\Cc} - 2 \ow{J_a\B}{J_a\Cc}\A -2
    \ow{J_a\Cc}\A{J_a\B}) \bigr).
  \end{multline*}
  Symmetrising in the first two variables we get the following element of
  $S^2V^*\otimes V^*$:
  \begin{equation*}
    \tfrac18\bigl(\vo\A\B\Cc-\vo\Cc\A\B + \sum_{a=1}^3(\vo\A{J_a\B}{J_a\Cc}
    - \vo\A{J_a\Cc}{J_a\B}) \bigr).
  \end{equation*}
  Applying the projection $p$ we get $1/32$ times
  \begin{multline*}
    \vo\A\B\Cc-\vo\Cc\A\B
    +\sum _{a=1}^3\bigl( \vo\A{J_a\B}{J_a\Cc} \\
    - \vo\A{J_a\Cc}{J_a\B}  + \vo{J_a\A}{J_a\B}\Cc -  \vo{J_a\Cc}{J_a\A}\B\\
    +\sum _{b=1}^3 (\vo{J_b\A}{J_bJ_a\B}{J_a\Cc}-
    \vo{J_bJ_a\Cc}{J_b\A}{J_a\B})\bigr).
  \end{multline*}
  Taking $\beta$ and $\gamma$ linearly independent over $\mathbb H$, which
  is possible for $\dim M \geqslant 8$, i.e., $n\geqslant2$, and examining
  the coefficient of $\cdot \otimes \gamma$, one sees that this element is
  non-zero.  Thus $\QK5\cong [KH]$ and hence $\QK4\cong [S^3EH]$.
\end{proof}

\begin{corollary}
  Only the following inclusions hold between classes of homogeneous
  quaternionic K\"ahler structures and homogeneous Riemannian structures:
  \begin{enumerate}
  \item $\QK5 \subset \mathcal T_2$,
  \item $\QK{2+4+5} \subset \mathcal T_{2+3}$.
  \end{enumerate}
  In particular, a naturally reductive homogeneous quaternionic K\"ahler
  structure is symmetric.
\end{corollary}

\begin{proof}
  Use the above descriptions of the modules $\QK i$ and Table~I, page~41
  in~\cite{Tricerri-V:homogeneous}.
\end{proof}

\section{Geometric results}
\label{sec:geo}

Let $(M, g, \upsilon)$ be a connected, simply-connected and complete
quaternionic K\"ahler manifold of dimension~$4n$.  Then each tangent space
$T_pM$, $p\in M$, with $(g, J_1, J_2, J_3)_p$ is a quaternion-Hermitian
vector space.  One has the standard representation of $\SP(n)\SP(1)$
on~$T_pM$ and hence it is possible to define and decompose the vector space
$\V_p \subset \mathcal{T}(T_pM)$ of pointwise homogeneous quaternionic
K\"ahler structures as in the previous section.  This decomposition depends
only on~$\upsilon_p$ and not on the chosen bases $(J_1,J_2,J_3)_p$, so the
irreducible summands~$(\QK i)_p$ give well-defined bundles~$\QK i$
over~$M$.

Suppose that $M$ admits a non-vanishing homogeneous quaternionic K\"ahler
structure~$S$.  Then, by Theorem~\ref{thqKhm}, $M$~is homogeneous.  Hence,
if $S_p$ belongs to a given invariant subspace of~$\V_p$, at $p \in M$,
then $S_q$ belongs to the similar invariant subspace of~$\V_q$ at any
other~$q \in M$ and is a section of the corresponding vector bundle.

\subsection{The class $\QK{1+2+3}$}

The purpose of this section is to prove one implication of
Theorem~\ref{thm:main}, namely:

\begin{theorem}
  \label{thm:linear}
  Suppose $M$ is a connected quaternionic K\"ahler manifold of dimension
  $4n\geqslant8$ admitting a non-vanishing homogeneous quaternionic
  structure $S\in\QK{1+2+3}$.  Then $S$ belongs to $\QK3$ and $M$~is
  locally isometric to the quaternionic hyperbolic space~$\HH(n)$.
\end{theorem}

The question of existence of $\QK3$-structures on~$\HH(n)$ will not be
addressed until~\S\ref{sec:qk3}.

\begin{remark}
  Computing dimensions one finds
  \begin{gather*}
    \dim[EH] = 4n, \quad  \dim[S^3E H] = 8n,  \\
    \dim[S^3EH] = \tfrac43n(n+1)(2n+1), \quad \dim[KH] =
    \tfrac{16}3n(n^2-1),
  \end{gather*}
  so $\QK1$, $\QK2$ and $\QK3$ are the modules whose dimensions grow
  linearly with $\dim M$ in $T^*\otimes (\sP(1) \oplus \sP(n))$.
  It is thus plausible that $\QK{1+2+3}$ corresponds to spaces of constant
  negative quaternionic curvature, since these are scarce in all
  homogeneous quaternionic K\"ahler spaces.  This phenomenon is similar to
  the Riemannian \cite{Tricerri-V:homogeneous} and K\"ahler
  \cite{Abbena-G:aH,Gadea-MA-MM:complex} cases.
\end{remark}

The proof of Theorem~\ref{thm:linear} will be divided in to a number of
steps.  For a couple of these we give alternate arguments: one
representation-theoretic, the other tensorial.  The former are often
shorter and more transparent; however, certain details of the tensorial
calculations are needed in the derivation of later results.

\begin{lemma}
  \label{lem:loc-symm}
  Suppose $M$ is a connected quaternionic K\"ahler manifold of dimension
  $4n \geqslant8$ with a non-vanishing homogeneous quaternionic K\"ahler
  structure $S \in \QK{1+2}$.  Then $M$ is locally symmetric.
\end{lemma}

\begin{proof}[First proof]
  The curvature~$R$ of any $4n$-dimensional quaternionic K\"ahler
  manifold~$M$ is given \cite{Alekseevsky:exceptional} by
  \begin{equation*}
    R = \nu_q R_0 + R_1,
  \end{equation*}
  with $R_0$ equal to the curvature tensor of the quaternionic projective
  space~$\HP(n)$, $R_1 \in [S^4E]$ and $\nu_q = \nu/4 =
  \mathbf{s}/16n(n+2)$, one quarter of the reduced scalar curvature.  If the
  structure is homogeneous, then
  \begin{equation*}
    0 = \Nt R = \nu_q \Nt R_0 + \Nt R_1  = \LC R_1 - SR_1
  \end{equation*}
  for the following reasons.  The tensor~$R_0$ is an
  $\SP(n)\SP(1)$-invariant algebraic curvature tensor built from the metric
  and local quaternionic structures in such a way that $\LC R_0=0$.  Also
  $S$~is an element of $\V=T^*M\otimes (\sP(1) \oplus \sP(n))$ and $S_X$
  acts via the differential of the $\SP(n)\SP(1)$-action, so $S_XR_0=0$.

  We may further decompose $S = S_H + S_E$, where $S_H \in T^*M\otimes
  \sP(1)$ and $S_E \in T^*M\otimes \sP(n)$.  As $R_1\in [S^4E]$, we have
  that $S_HR_1=0$.  Thus $\LC R = \LC R_1 = S_ER_1$.  We conclude that if
  $S_E=0$, i.e., if $S$~is of type~$\QK{1+2}$, then $\LC R=0$ and $g$~is
  locally symmetric.
\end{proof}

\begin{proof}[Second proof]
  Any quaternionic K\"ahler manifold of dimension $4n \geqslant8$ is
  Einstein~\cite{Alekseevsky:exceptional,Berger:CR}.  Moreover, one has
  \begin{multline}
    \label{jjqk}
    R_{XYJ_1ZW} + R_{XYZJ_1W} \\
    = \tfrac1{n+2} \bigl( \mathbf{r}(J_2X,Y)g(J_3Z,W)
    - \mathbf{r}(J_3X,Y)g(J_2Z,W) \bigr),
  \end{multline}
  etc., which is proved in \cite[(2.13)]{Ishihara:qK} under a different
  curvature convention, and which may be found with a misprint
  in~\cite[p.~404]{Besse:Einstein}.

  As $S\in\QK{1+2}$ we have locally that $S_{XYZ}=\theta^a(X)g(J_aY,Z)$,
  with $\theta^a\in T^*M$.  The second Ambrose-Singer equation $\Nt R=0$
  of~\eqref{eq:AS} is
  \begin{multline}
    \label{asiiud}
    (\LC_X R)_{YZWU} \\= - \theta^a(X) \bigl( R_{J_aYZWU} + R_{YJ_aZWU} +
    R_{YZJ_aWU} + R_{YZWJ_aU} \bigr).
  \end{multline}
  But~\eqref{jjqk} implies that the right-hand side vanishes.
\end{proof}

\begin{lemma}
  Let $M$ be a quaternionic K\"ahler manifold of dimension $4n\geqslant8$.
  Suppose $S$ is a homogeneous quaternionic K\"ahler structure with $S \in
  \QK{1+2+3}$ and with non-zero projection to $\QK3$. Then $M$ has constant
  quaternionic curvature.
\end{lemma}

\begin{proof}[First proof]
  We saw in the first proof of Lemma~\ref{lem:loc-symm} that $\LC R = S_E
  R_1$.  The assumption that the projection of~$S$ to~$\QK3$ is non-zero
  implies that~$S_E\in[EH]$ and is non-zero.  By the differential Bianchi
  identity one has that $\LC R \in [S^5EH]$, see~\cite[proof of
  Th.~2.6]{Poon-Salamon:8}, thus
  \begin{equation}
    \label{eq:SER1}
    S_E\otimes R_1 \in [EH] \otimes [S^4E] \cong [S^5EH] + [S^3EH] + [V^{(31)}H],
  \end{equation}
  where $V^{(31)}$ is irreducible.

  The map $S_E \otimes R_1 \mapsto S_ER_1$ is the composition
  \begin{equation*}
    \phi\colon [EH] \otimes [S^4E] \hookrightarrow [EH] \otimes
    [S^2E] \otimes [S^4E] \longrightarrow [EH] \otimes [S^4E],
  \end{equation*}
  where the first map is given by the inclusion of $[EH]$
  in~$[EH]\otimes[S^2E]$ as the module~$\QK3$ and the second map is given
  by the action of $[S^2E]\cong\sP(n)$ on~$[S^4E]$.  This composition is
  linear and $\SP(n)\SP(1)$-equivariant.

  Write $\{h,\tilde h\}$ for a complex orthogonal basis of $H$.  Consider
  the elements $\alpha = e_1h\otimes {\tilde e_1}^4 $ and $\beta = e_1h
  \otimes e_1 \vee \tilde e_1 \vee e_2^2$ of $EH\otimes S^4E$, where
  $\tilde e=je$, and $e_1,\tilde e_1,e_2$ are linearly independent.  Then
  $\phi(\alpha)$ has non-zero components in $S^5EH$ and $S^EH$, whereas
  $\phi(\beta)$ lies in neither of these modules.  By Schur's Lemma, we
  conclude that $\phi$ is an isomorphism on each component of the
  decomposition~\eqref{eq:SER1}.

  As $\nabla R=S_ER_1\in[S^5EH]$, we have $S_E\otimes R_1\in[S^5H]$ too.
  Now $S_E = eh + \tilde e\tilde h$, with $\tilde e=je$, since $S_E$~is a
  real element.  For $S_E\otimes R_1$ to be totally symmetric in the $e$'s
  we must have that $R_1 \in S^4\{e,\tilde e\}$.  But $h$ and $\tilde h$
  are linearly independent so $e \otimes R_1$ and $\tilde e \otimes R_1$
  are each in $S^5E$; the first implies that $R_1 = a e^5$, the second that
  $R_1 = b {\tilde e}^5$.  As $e$ and $\tilde e$ are linearly independent,
  we conclude that $R_1=0$.  Thus $R=\nu_qR_0$ and our space has constant
  quaternionic curvature.
\end{proof}

\begin{proof}[Second proof]
  For $S \in \QK{1+2+3}$, we have locally
  \begin{equation}
    \label{eq:QK123}
    \begin{split}
      S_XY &= g(X,Y) \xi- g(\xi,Y) X \\
      &\qquad +\sum _{a=1}^3\bigl( g(\xi,J_a Y) J_a X - g(X,J_a Y) J_a \xi
      + g(X, \zeta ^a) J_aY\bigr),
    \end{split}
  \end{equation}
  with $\xi$ and $\zeta^a$ vector fields metrically dual to the one-forms
  $\theta$ and $\theta^a$ of Theorem~\ref{tchqKs}.  Our assumptions imply
  in addition that $\xi\ne0$.

  We compute first, $\Nt\xi$ and $\Nt\zeta ^a$.  The third Ambrose-Singer
  equation of~\eqref{eq:AS} can be written as $\Nt_Z(S_X Y)= S_{\Nt_ZX}Y +
  S_X\Nt_ZY$.  Taking the covariant derivative of~\eqref{eq:QK123} with respect
  to~$Z$ and using equations~\eqref{nabj}, we get
  \begin{equation}
    \label{eq:Nt-Z}
    \begin{split}
      0 &= g(X,Y)\Nt_Z\xi - g(\Nt_Z\xi,Y)X\\
      &\qquad + \sum _{a=1}^3 \bigl( g(\Nt_Z\xi, J_a Y) J_aX
      + g(X,J_a Y) J_a\Nt_Z\xi \bigr) \\
      &\qquad + \Cyclic_{123} g\bigl(X, \Nt_Z\zeta^1 -
      \tilde\tau^3(Z)\zeta^2 + \tilde\tau^2(Z)\zeta^3\bigr) J_1Y,
    \end{split}
  \end{equation}
  for all $X,Y,Z$.  We may perform two operations on this equation.  Either
  take the inner product with~$X$ or put $X=Y$.  In each case, now take the sum
  over $X$ in an orthonormal basis.  This gives the two relations
  $(4n-4)\Nt_Z\xi=\pm W$, for a certain vector field $W$, and we thus have
  $\Nt_Z\xi=0$, or, equivalently,
  \begin{equation*}
    \LC_Z\xi=S_Z\xi.
  \end{equation*}
  It follows that $g(\xi,\xi) $ is a constant function, as $Xg(\xi,\xi) =
  2g(\Nt_X\xi,\xi) = 0$ for all $X\in\mathfrak{X}(M)$.
  Equation~\eqref{eq:Nt-Z} now implies
  \begin{equation}
    \label{eq:nabz}
    \Nt_Z\zeta^1 = \tilde\tau^3(Z)\zeta^2 - \tilde\tau^2(Z)\zeta^3,
    \quad\text{etc.}
  \end{equation}

  The second Ambrose-Singer equation of~\eqref{eq:AS} can be written as
  \begin{equation}
    \label{asiiudt}
    \begin{split}
      &(\LC_XR)_{YZWU}\\
      &\qquad = -R_{S_XYZWU} - R_{YS_XZWU} - R_{YZS_XWU} - R_{YZWS_XU}.
    \end{split}
  \end{equation}
  Substituting \eqref{eq:QK123} in~\eqref{asiiudt}, one sees that the terms
  containing~$\zeta^a$ are expressed as the summands with four $R$'s in the
  right-hand side of~\eqref{asiiud}, but we have seen that each such
  summand vanishes.

  Taking then the cyclic sum with respect to $X,Y,Z$ in~\eqref{asiiudt}, we
  get after a quite long calculation, using the two Bianchi
  identities and the relation~\eqref{jjqk}, that
  \begin{multline}
    \label{35a}
    R_{XY}\xi = \nu_q \bigl\{ g( X,\xi) Y-g( Y,\xi) X \\ + \sum
    _{a=1}^3(g( J_a X,\xi) J_a Y -g(J_a Y, \xi) J_a X - 2g( J_a X,Y) J_a
    \xi ) \bigr\},
  \end{multline}
  which is the expression of the curvature tensor~$R(X,Y)Z$, for $Z = \xi$,
  of a quaternionic K\"ahler manifold with constant quaternionic curvature
  $\nu=4\nu_q$ \cite{Alekseevsky:exceptional,Ishihara:qK}.  We only need to
  prove that the expression similar to~\eqref{35a} is true for $R_{XY}Z$,
  with $Z$ arbitrary, instead of the particular $\xi$.

  For this, we apply again the second Bianchi identity to the second
  Ambrose-Singer equation $\Nt R=0$, so
  \begin{equation}
    \label{44}
    0= \Cyclic_{XYZ} (R_{S_XYZWU}+R_{YS_XZWU}+R_{YZS_XWU} +R_{YZWS_XU}).
  \end{equation}
  For the sake of simplicity, we write $\Theta^1(X,Y,Z,W)$ for the
  right-hand side of~\eqref{jjqk}.  Expanding the terms $S_X$
  in~\eqref{44}, on account of formul\ae\ \eqref{eq:QK123} and~\eqref{jjqk}, we
  obtain that
  \begin{equation}
    \label{55}
    \begin{split}
      0 &
      = \Cyclic_{XYZ}
      \bigl\{-2g( X,\xi ) R_{ZYWV} -2g(X,J_a Y ) R_{J_a \xi ZWV} +g( X,W)
      R_{YZ\xi U} \\
      &\quad\qquad +g( X,U) R_{YZW\xi} -g( X,J_a W) R_{YZJ_a
      \xi U}- g(X,J_a U) R_{YZWJ_a \xi}\bigr\} \\
      & \quad + \Cyclic_{XYZ} \bigl\{  g(\xi,J_a X)\Theta ^a
      (W,U,Z,Y)+g(J_a W,\xi)\Theta ^a (X,Y,Z,U) \\
      &\quad\qquad +g( J_a U,\xi)\Theta ^a(X,Y,W,Z) -2g( X,J_a Y)\Theta ^a
      (W,U,\xi, Z)\\
      &\quad\qquad -g( X,J_a W) \Theta ^a (Y,Z,\xi, U)-g( X,J_a U)\Theta ^a
      (Y,Z,W,\xi)\bigr\},
    \end{split}
  \end{equation}
  where, as we saw, the terms in $\zeta ^a$ in \eqref{eq:QK123} do not actually
  contribute.

  Again from \eqref{jjqk}, after some computations, the second cyclic sum
  above can be written as
  \begin{equation*}
    \nu\,\Cyclic_{123}\bigl(  \omega_2 (W,U)(J_3 \xi)^\flat -\omega_3
    (W,U)(J_2 \xi)^\flat \bigr) \wedge\omega_1.
  \end{equation*}

  We now work out the first cyclic sum in~\eqref{55}.  First, we write this
  as
  \begin{multline*}
    2\theta\wedge R_{WU} - 2\sum_{a=1}^3\omega_a \wedge R_{WU\xi J_a
    (\cdot)} +W^\flat\wedge R_{\xi U} -U^\flat\wedge
    R_{\xi W} \\
    + \sum _{a=1}^3\bigl((J_a W)^\flat \wedge R_{\xi J_aU} - (J_a
    U)^\flat \wedge R_{\xi J_aW}\bigr).
  \end{multline*}
  Now making use of formula~\eqref{35a}, after some simplifications we
  obtain that this can be written as
  \begin{multline*}
    \theta \wedge \bigl\{ 2R_{WU} - \tfrac12 \nu\,\bigl[ W^{\flat }\wedge
    U^\flat + \sum _{a=1}^3\bigl((J_a W)^\flat\wedge(J_a U)^\flat + 2
    \omega_a
    (W,U) \omega_a \bigr) \bigr] \bigr\} \\
    + \nu \Cyclic_{123}\bigl\{ \omega_3 (W,U) (J_2 \xi)^\flat -\omega_2
    (W,U)(J_3 \xi)^\flat\bigr\} \wedge\omega_1.
  \end{multline*}
  Hence, expression~\eqref{55} can be written as
  \begin{equation*}
    0=\theta\wedge\bigl\{ R_{WU} - \nu_q\bigl[ W^{\flat }\wedge
    U^\flat + \sum _{a=1}^3\bigl((J_a W)^\flat \wedge(J_a U)^\flat + 2
    \omega_a (W,U)\omega_a \bigr) \bigr] \bigr\}.
  \end{equation*}
  Contracting with $\xi$, we obtain that
  \begin{equation*}
    \begin{split}
      0 & =\Bigl\{ R_{WU} - \nu_q \bigl[ W^\flat\wedge U^\flat +
      \sum_{a=1}^3\bigl((J_a W)^\flat\wedge (J_a U)^\flat + 2 \omega_a
      (W,U)\omega_a
      \bigr)\bigr] \Bigr\} \\
      &\quad - \theta \wedge \biggl\{ R_{WU\xi} - \nu_q
      \Bigl[ g( W,\xi) U^\flat - g( U,\xi) W^\flat
      + \sum _{a=1}^3\bigl(g( J_a W,\xi)(J_a U)^\flat\\
      &\hspace{10em} - g( J_a U,\xi) (J_a W)^\flat+ 2 \omega_a
      (W,U)(J_a\xi)^\flat\bigr)\Bigr]\biggr\}.
    \end{split}
  \end{equation*}
  By formula~\eqref{jjqk}, the second curly bracket vanishes, so
  \begin{equation*}
    R_{WU} = \nu_q \{W^\flat
    \wedge U^\flat +\sum _{a=1}^3\bigl( (J_a W)^\flat
    \wedge(J_a U)^\flat +2 \omega_a (W,U) \omega_a \bigr)\},
  \end{equation*}
  i.e., $(M,g,\upsilon)$ is a space of constant quaternionic
  curvature~$\nu=4\nu_q$.
\end{proof}

\begin{proposition}
  \label{prop:QK3}
  If $S \in \QK{1+2+3}$ and has non-zero projection to~$\QK3$, then the
  manifold is locally isometric to the quaternionic hyperbolic space and
  $S$ belongs to~$\QK3$.
\end{proposition}

\begin{proof}
  By hypothesis, the tensor~$S$ is given by~\eqref{eq:QK123} with $\xi\ne0$, from
  which $\Nt\xi=0$ and equation~\eqref{35a} were derived, however we do not
  yet know the value of~$\nu_q$.  On the other hand, using $\LC=\Nt+S$,
  $\Nt\xi=0$ and equations~\eqref{naji} and~\eqref{nabj}, we have
  \begin{equation}
    \label{tresquince}
    \begin{split}
      \LC_XJ_1\xi
      &= \tilde\tau^3(X)J_2\xi - \tilde\tau^2(X)J_3\xi +g( X,J_1\xi )
      \xi\\
      &\quad +\sum _{b=1}^3 \bigl(g( \xi ,J_bJ_1\xi ) J_bX
      -g( X, J_bJ_1\xi ) J_b\xi +g( X,\zeta ^b ) J_bJ_1\xi\bigr),
    \end{split}
  \end{equation}
  which will be used in conjunction with
  \begin{equation}
    \label{seis}
    \begin{split}
      R_{XY}\xi
      &= - (\LC_X(\LC\xi))_Y + (\LC_Y(\LC\xi))_X
      = - (\LC_X(S\xi))_Y + (\LC_Y(S\xi))_X\\
      &= - g(Y,\LC_X\xi)\xi - g(Y,\xi)\LC_X\xi + g(X,\LC_Y\xi)\xi +
      g(X,\xi)\LC_Y\xi\\
      & \quad +\sum _{a=1}^3 \bigl(g(Y,\LC_XJ_a\xi)J_a\xi +
      g(Y,J_a\xi)\LC_XJ_a\xi \\ 
      & \qquad\qquad
      - g(X,\LC_YJ_a\xi) J_a\xi - g( X,J_a\xi)\LC_YJ_a\xi \\
      & \qquad\qquad - g(Y,\LC_X\zeta^a)J_a\xi - g(Y,\zeta^a)\LC_XJ_a\xi \\
      & \qquad\qquad + g(X,\LC_Y\zeta^a)J_a\xi + g(X,\zeta^a )\LC_YJ_a\xi\bigr).
    \end{split}
  \end{equation}
  This will be examined in three stages to find the sign of the scalar
  curvature and to show that each $\zeta^a$ is zero.

  \textbf{Step~1}: Take $X,Y\in (\mathbb H\xi )^{\perp }$.  Then
  \eqref{35a} says $g(R_{XY}\xi,X)=0$.  However, from~\eqref{seis} we have
  \begin{equation*}
    g(R_{XY}\xi,X) = - g(Y,\zeta^a)g(\LC_XJ_a\xi,X) +
    g(X,\zeta^a)g(\LC_YJ_a\xi,X).
  \end{equation*}
  Now taking $Y=J_aX$ and using~\eqref{tresquince}, we find that
  $g(R_{XJ_aX}\xi,X) = g(X,\zeta^a)g(\xi,\xi)g(X,X)$.  As $X \in (\mathbb
  H\xi)^\perp$ is arbitrary, one deduces that
  \begin{equation}
    \label{ss1}
    \zeta^a \in \mathbb H\xi.
  \end{equation}

  \textbf{Step~2}: Take $X=\xi $ and $Y\in (\mathbb H\xi)^{\perp }$
  in~\eqref{seis}.  As $g(\xi,\xi)$~is constant, we have that
  $g(\xi,\LC_Y\xi)=0$.  Moreover, using~\eqref{naji} we find
  \begin{gather*}
    g(Y,\LC_{\xi}\xi) = 0,\quad
    g(Y,\LC_{\xi}\zeta^a) = 0,\\
    g(Y,\LC_\xi J_1\xi) = g(Y,J_1\LC_\xi\xi) + g(Y,\tau^3(\xi)J_2\xi) -
    g(Y,\tau^2(\xi)J_3\xi) = 0,\\
    g(\xi,\LC_YJ_1\xi) = g(\xi,J_1\LC_Y\xi) + g(\xi,\tau^3(Y)J_2\xi) -
    g(\xi,\tau^2(Y)J_3\xi) = 0,
  \end{gather*}
  etc., which leads to
  \begin{equation*}
    R_{\xi Y}\xi = - g(\xi,\xi)^2Y-g(\xi,\zeta^a)g(\xi,\xi)J_aY.
  \end{equation*}
  Comparing with~\eqref{35a} which says $R_{\xi Y}\xi = \nu_q g(\xi,\xi)Y$,
  we obtain $\nu_q=-g(\xi,\xi)$ and
  \begin{equation}
    \label{ss2}
    g(\xi,\zeta^a) = 0.
  \end{equation}

  \textbf{Step~3}: Take $X,Y\in (\mathbb H\xi )^{\perp }$  again and
  use~\eqref{seis} to get
  \begin{equation*}
    \begin{split}
      R_{XY}\xi
      &= - g(Y,\LC_X\xi)\xi + g(X,\LC_Y\xi)\xi\\
      &\quad+\sum_{a=1}^3 \bigl( g(Y,\LC_XJ_a\xi)J_a\xi-g(X,\LC_YJ_a\xi)J_a\xi\\
      &\qquad\qquad -g(Y,\LC_X\zeta^a)J_a\xi+g(X,\LC_Y\zeta^a)J_a\xi \bigr).
    \end{split}
  \end{equation*}
  Making use of \eqref{eq:nabz}, the expression for~$S$ and \eqref{ss1}, we
  obtain after some calculations that
  \begin{equation*}
    \begin{split}
      R_{XY}\xi
      &= - 2 \sum_{a=1}^3 g(Y,J_aX)g(\xi,\xi)J_a\xi\\
      &\qquad - \sum_{a,b=1}^3
      \bigl(g(Y,J_aX)g(\xi,J_a\zeta^b)
      - g(X,J_aY)g(\xi,J_a\zeta^b)\bigr)J_b\xi.
    \end{split}
  \end{equation*}
  However, formula~\eqref{35a} says $R_{XY}\xi = -2\sum_{a=1}^3 g(\xi,\xi)
  g(J_aX,Y) J_a\xi$.  So we have
  \begin{equation*}
    \sum_{a=1}^3 g(Y,J_aX)g(\xi,J_a\zeta^b) - g(X,J_aY)g(\xi,J_a\zeta^b) =
    0,\quad\text{for $b=1,2,3$}.
  \end{equation*}
  Taking $Y=J_cX$ we conclude that $-2g(X,X)g(\xi,J_c\zeta^b) = 0$, for all
  $b,c$.  Together with \eqref{ss1} and~\eqref{ss2} this gives $\zeta^a=0$,
  for each $a$, and hence $S\in\QK3$.  Consequently, our manifold is
  locally isometric to the hyperbolic space of constant quaternionic
  curvature $-4g(\xi,\xi)$, as claimed.
\end{proof}

\subsection{Non-existence of $\QK{1+2}$}

To complete the proof of Theorem~\ref{thm:linear} we need to show that
structures of type $\QK{1+2}$ do not occur on manifolds of dimension $8$ or
more.  For these structures, the tensor $S$~lies in $[EH] \otimes [S^2H]$.

\begin{lemma}
  If a $(1,2)$ tensor $S$ satisfies the conditions for the tensors in
  $\QK{1+2}$ and also $\Nt S = 0$, then $S$ defines a homogeneous
  structure.
\end{lemma}

\begin{proof}
  We need to show that $\Nt\Rt = 0$.  Since $\Nt S = 0$, we have that
  \begin{equation}
    \label{eq:Rt-R}
    \Rt = R - R^S = \nu_q R_0 + R_1 - R^S,
  \end{equation}
  where $R_0$ is the curvature tensor of~$\HP(n)$, $R_1 \in [S^4E]$ and
  \begin{equation}
    \label{eq:RS}
    R^S_{XY}Z = S_Y(S_XZ) - S_{S_YX}Z - (S_X(S_YZ)-S_{S_XY}Z).
  \end{equation}
  By Lemma~\ref{lem:loc-symm}, $M$ is a locally symmetric space, so $\LC
  R_1=0$.  This implies $\Nt_X R_1 = - S_XR_1 =0$, since $[S^2H]\cong
  \sP(1)$ acts trivially on~$[S^4E]$.  However, $\Nt$~is an
  $\SP(n)\SP(1)$-connection, so $\Nt_X R_0=0$.  We are left with $\Nt\Rt =
  -\Nt R^S$, which is zero once $\Nt S=0$.
\end{proof}

Since $S\in [EH]\otimes[S^2H]$ we may write $S=(X^a)^\flat\otimes J_a$ for
some vector fields $X^a$, which means $S_X = g(X,X^a)J_a$.  The equation
$\Nt_X S =0$ is equivalent to $\LC_X S =S_X.S$.  Now
\begin{equation*}
  \begin{split}
    (S_X.S)_Y Z
    &= S_XS_YZ - S_{S_XY}Z- S_YS_XZ\\
    &= \Cyclic_{123}\bigl\{2 g(X,X^2) g(Y,X^3) - 2 g(X,X^3) g(Y,X^2)\\
    &\qquad\qquad -g(X,X^b) g(J_bY,X^1) \bigr\}J_1Z.
  \end{split}
\end{equation*}
Thus $\Nt S=0$ is equivalent to
\begin{multline}
  \label{eq:Nt-S}
  \LC X^1 = (2(X^2)^\flat-\tau^2)\otimes X^3
  - (2(X^3)^\flat-\tau^3)\otimes X^2 \\
  + (X^a)^\flat \otimes J_aX^1,
\end{multline}
etc., where $\tau^a$ are given by~\eqref{naji}.

\begin{lemma}
  If $M$ admits a non-vanishing homogeneous structure in the class
  $\QK{1+2}$, then $\dim M\leqslant 12$.  Moreover $M$ is not of type
  $\QK1$.
\end{lemma}

\begin{proof}
  Equation~\eqref{eq:Nt-S} implies that
  \begin{equation*}
    \LC_X X^1 \in \Span\{X^2,X^3,J_1X^1,J_2X^1,J_3X^1\},\quad\text{etc.}
  \end{equation*}
  So the distribution given by the quaternionic span $\mathcal D=\mathbb
  H\left\{X^1,X^2,X^3\right\}$ is parallel and hence holonomy
  invariant~\cite[Prop.\ 10.21]{Besse:Einstein}.  The irreducibility of $M$
  implies that this distribution must be the whole tangent space.  This
  implies that $\dim_{\mathbb R}M\leqslant12$.  If $M$ is of type $\QK1$,
  then $X^a=J_aX_0$ for a fixed non-zero vector field~$X_0$.  But this
  implies that $\mathcal D=TM$ is four-dimensional, which we have
  specifically excluded.
\end{proof}

Note that we may now assume that at least two of the $X^a$'s are linearly
independent over~$\mathbb H$.

Let us now compute the Riemannian curvature.  Regarding $R$ as the
alternation of $-\nabla\nabla$ we compute
\begin{equation*}
  \begin{split}
    RX^1
    &= -\mathbf a\LC\LC X^1
    = -\mathbf a\LC(2X^2- \tau^2)X^3 + \mathbf a(2X^2-\tau^2)\nabla X^3 +
    \dotsb\\
    &= (2\kappa^3-\Omega^3)X^2 - (2\kappa^2-\Omega^2)X^3 - \kappa^1J_1X^1
    - \kappa^2J_2X^1 - \kappa^3J_3X^1,
  \end{split}
\end{equation*}
where
\begin{equation*}
  \Omega^a = d\tau^a + \tau^b\wedge\tau^c,\qquad
  \kappa^a = 2X^b\wedge X^c + X^d\wedge J_dX^a,
\end{equation*}
with $(a,b,c)$ a cyclic permutation of $(1,2,3)$.

As $M$ is quaternionic K\"ahler, the curvature splits as
\begin{equation*}
  R = R_E + R_H,
\end{equation*}
with $R_E \in [S^4E] + \mathbb R \subset [S^2(S^2E)]$ and $R_H \in
\mathbb R \subset [S^2(S^2H)]$.  The part with values in $[S^2E]$ is
given by
\begin{equation*}
  R_E X = \tfrac14(R X - \sum_a J_aR(J_aX)) = RX - \tfrac 14\sum_a
  J_aR(J_a)X.
\end{equation*}
Now $R J_1 = - \Omega^3J_2 + \Omega^2J_3$, so
\begin{equation*}
  \begin{split}
    R_E X^1 &= R X^1 - \tfrac12\left(\Omega^1J_1 + \Omega^2J_2 +
      \Omega^3J_3 \right)X^1\\
    &= 2\Theta^3X^2 -2\Theta^2X^3 - \Theta^1J_1X^1 - \Theta^2J_2X^1 -
    \Theta^3J_3X^1,
  \end{split}
\end{equation*}
where $\Theta^a = \kappa^a - \tfrac12\Omega^a$.

Let us work with the above expression in two complementary ways.  Firstly,
we have that $(R_E)_{AB}$ is a skew-adjoint endomorphism of~$TM$.  Thus
\begin{equation*}
  0=\tfrac12\inp{R_EX^1}{X^1} = \alpha^3\Theta^3 - \alpha^2\Theta^2,
\end{equation*}
where
\begin{equation*}
  \alpha^a = \inp{X^b}{X^c}.
\end{equation*}
And similarly
\begin{equation*}
  0
  = \tfrac12\inp{R_EX^1}{X^2} + \tfrac12\inp{R_EX^2}{X^1}
  = \alpha^2\Theta^1 - \alpha^1\Theta^2 + \Theta^3(\delta^2-\delta^1),
\end{equation*}
where
\begin{equation*}
  \delta^a = \norm{X^a}^2.
\end{equation*}
In addition we know that $R_E$ commutes with each $J_a$, so
\begin{equation*}
  \begin{split}
    0 &= \tfrac12\inp{R_E(J_1X^1)}{X^2} + \tfrac12\inp{R_EX^2}{J_1X^1}\\
    &= -\tfrac12\inp{R_EX^1}{J_1X^2} + \tfrac12\inp{R_EX^2}{J_1X^1} \\
    &= \beta^{13}\Theta^1 + (\gamma^3-\gamma^2)\Theta^2 -
    \beta^{31}\Theta^3,
  \end{split}
\end{equation*}
where
\begin{equation*}
  \beta^{bc} = \inp{J_bX^b}{X^c}\quad \text{no sum},\qquad \gamma^a =
  \inp{X^a}{J_bX^c}.
\end{equation*}
Similarly,
\begin{equation*}
  0 = \tfrac12\inp{R_E{J_2X^1}}{X^2} + \tfrac12\inp{R_EX^2}{J_2X^1}
  = (\gamma^2-\gamma^3)\Theta^1 + \beta^{23}\Theta^2 - \beta^{12}\Theta^3
\end{equation*}
and
\begin{equation*}
  0 = \tfrac12\inp{R_E{J_3X^1}}{X^2} + \tfrac12\inp{R_EX^2}{J_3X^1}
  = (\beta^{21}-\beta^{31})\Theta^1 + (\beta^{12}-\beta^{32})\Theta^2.
\end{equation*}
We thus have the following system of equations
\begin{gather}
  \label{eq:a-T}
  \alpha^1\Theta^1 = \alpha^2\Theta^2 = \alpha^3\Theta^3,\\
  \label{eq:a-d}
  \alpha^b\Theta^a - \alpha^a\Theta^b + (\delta^b-\delta^a)\Theta^c = 0,\\
  \label{eq:b-g1}
  \beta^{ac}\Theta^a + (\gamma^c-\gamma^b)\Theta^b - \beta^{ca}\Theta^c = 0, \\
  \label{eq:b-g2}
  (\gamma^b-\gamma^c)\Theta^a + \beta^{bc}\Theta^b - \beta^{ab}\Theta^c = 0, \\
  \label{eq:b-g3}
  (\beta^{ba}-\beta^{ca})\Theta^a + (\beta^{ab}-\beta^{cb})\Theta^b = 0,
\end{gather}
for each cyclic permutation $(a,b,c)$ of $(1,2,3)$.

Now consider the Riemannian holonomy of~$M$.  We know that $M$ is locally
symmetric and quaternionic K\"ahler.  Thus the holonomy algebra $\hol$ is
that of a quaternionic symmetric space and splits as
\begin{equation*}
  \hol = \lie k \oplus \sP(1) \subset \sP(n) \oplus \sP(1),
\end{equation*}
with $\lie k$ non-Abelian.  However, the Lie algebra $\lie k$ is generated
by the coefficients of~$R_E$.  In particular, the linear span of
$\Theta^1$, $\Theta^2$ and $\Theta^3$ has to be at least two-dimensional.
Equations~\eqref{eq:a-T} and~\eqref{eq:a-d}, then imply that
\begin{equation*}
  \alpha^a=0,\qquad \delta^1=\delta^2=\delta^3,
\end{equation*}
so the $X^a$ are mutually orthogonal and of equal length.

\begin{lemma}
  \begin{equation*}
    \gamma^1=\gamma^2=\gamma^3,\qquad \beta^{bc}=0.
  \end{equation*}
\end{lemma}

\begin{proof}
  If $\dim_{\mathbb R}\Span\{\Theta^1,\Theta^2,\Theta^3\} = 3$, then
  this is direct from equations \eqref{eq:b-g1}, \eqref{eq:b-g2}
  and~\eqref{eq:b-g3}.

  If $\Theta^a$ are not linearly independent, we may without loss of
  generality assume
  \begin{equation*}
    \Theta^1 = x\Theta^2 + y\Theta^3
  \end{equation*}
  with $\Theta^2$ and $\Theta^3$ linearly independent.  If $x$
  and~$y$ are both zero, then the result follows easily.  Assume
  therefore that $x\ne0$.  We get from~\eqref{eq:b-g3},
  \begin{equation*}
    \beta^{12}=\beta^{32},\quad \beta^{23}=\beta^{13},\quad
    \beta^{21}=\beta^{31}.
  \end{equation*}
  Comparing the coefficient of $\Theta^2$ in the first equation
  of~\eqref{eq:b-g1} with the first equation of~\eqref{eq:b-g2} gives
  \begin{equation*}
    \beta^{23}x = \gamma^2-\gamma^3,\qquad (\gamma^2-\gamma^3)x = -\beta^{23},
  \end{equation*}
  from which one concludes $\beta^{23}=0$ and $\gamma^2=\gamma^3$.
  The result follows.
\end{proof}

Write $\gamma=\gamma^a$ and $\delta=\delta^a$ which are now both
independent of the index~$a$.  We find that $X^1$ is orthogonal to
\begin{equation*}
  J_1X^1, J_2X^1, J_3X^1, X^2, J_1X^2, J_2X^2, X^3, J_1X^3, J_3X^3.
\end{equation*}
Note using equation~\eqref{eq:Nt-S} with this information shows that
$\delta=\norm{X^1}^2$ is constant.  We may write
\begin{equation}
  \label{eq:W1}
  X^1 = aJ_3X^2 + bJ_2X^3 + W^1,
\end{equation}
with $W^1$ orthogonal to the quaternionic span of $X^2$ and $X^3$.  Now
\begin{equation}
  \label{eq:g1}
  \gamma = \inp{X^1}{J_2X^3} = a\inp{J_3X^2}{J_2X^3} +
  b\inp{J_2X^3}{J_2X^3}
  = -a\gamma + b\delta.
\end{equation}
Also
\begin{equation}
  \label{eq:g2}
  \gamma = \inp{X^2}{J_3X^1} = -a\inp{X^2}{X^2} - b\inp{X^2}{J_1X^3} =
  -a\delta + b\gamma
\end{equation}
and
\begin{equation}
  \label{eq:delta}
  \delta = \norm*{X^1}^2 = (a^2+b^2)\delta - 2ab\gamma + \norm*{W^1}^2.
\end{equation}

Now the other way of looking at $R_E$ is to note that
$\inp{(R_E)_{\cdot,\cdot}X}Y$ is in $S^2E$ for all $X$ and $Y$.  This says
that these two-forms are of type $(1,1)$ for each~$J_a$.  Grouping $\Theta$
terms, we have
\begin{equation*}
  R_E X^1 = - \Theta^1 J_1X^1 - (2X^3+J_2X^1) \Theta^2 + (2X^2-J_3X^1) \Theta^3.
\end{equation*}
Now by the above analysis the vectors $J_1X^1$, $2X^3+J_2X^1$,
$2X^2-J_3X^1$ are mutually orthogonal and all non-zero.  Thus each
$\Theta^a$ lies in $S^2E$.  The condition that $\Theta^1 =
\kappa^1-\tfrac12\Omega^1$ is type $(1,1)_1$, i.e., type $(1,1)$ with
respect to $J_1$, tells us that
\begin{equation*}
  2X^2\wedge X^3 + X^2\wedge J_2X^1 + X^3\wedge J_3X^1
\end{equation*}
is of type $(1,1)_1$, since $\Omega^1$ is proportional to the
K\"ahler form $\omega_1$ for~$J_1$, and the term $X^1\wedge
J_1X^1$ is already of type $(1,1)_1$.  Squaring, we have that
\begin{equation*}
  J_2X^1 \wedge J_3X^1 \wedge X^2 \wedge X^3
\end{equation*}
is a $(2,2)_1$-form, which implies that $\Span\{ J_2X^1, J_3X^1, X^2,
X^3 \}$ is $J_1$-invariant.  In particular, $J_1X^2$ is a linear
combination of $J_2X^1$, $J_3X^1$ and~$X^3$.  We conclude that
the quaternionic span of $X^1$, $X^2$ and $X^3$ has $8$
real dimensions.

In equation~\eqref{eq:W1}, we find that $W^1=0$.  Note that
\begin{equation}
  \label{eq:strict}
  \abs{\gamma} = \abs{\inp*{J_1X^1}{J_2X^2}} < \delta,
\end{equation}
since the three vectors $J_1X^1$, $J_2X^2$ and $J_3X^3$ can not all be
proportional.  Adding equations~\eqref{eq:g1} and~\eqref{eq:g2}, we get
$(a+b)\gamma=(a+b)\delta$, so $a=-b$ by~\eqref{eq:strict}.  Now
by~\eqref{eq:delta}, $2a^2(\gamma+\delta) = \delta$ whereas \eqref{eq:g1}
gives $a(\gamma+\delta) = - \gamma$.  Thus $2a\gamma = -\delta$ and
$2a^2-a=1$.  The latter has solutions $a=-1/2$ and $a=1$.  However,
$a=-1/2$ is impossible by~\eqref{eq:strict}, so $a=1$, $b=-1$,
$\gamma=-\delta/2$.  In particular,
\begin{equation*}
  X^1 = J_3X^2 - J_2X^3,
\end{equation*}
which says the structure is of type~$\QK2$.

However, we will see that these structures do not arise.  Assume that the
constant $\delta$ is non-zero.  Rescaling the geometry by a homothety
we may assume that $\delta=1$, $\gamma=-1/2$.  Put
\begin{equation*}
  A = J_1X^1,\qquad B = \tfrac1{\sqrt3}(J_2X^2-J_3X^3).
\end{equation*}
This is an orthonormal quaternionic basis for $TM$.  We have
\begin{equation*}
  J_2X^2 = -\tfrac12 A + \tfrac{\sqrt3}2 B,\qquad
  J_3X^3 = -\tfrac12 A - \tfrac{\sqrt3}2 B.
\end{equation*}
For $\Theta^1$ to be of type $(1,1)_2$, we must have that
$\kappa^1-J_2\kappa^1 = \Omega^1$.  But
\begin{equation*}
  \begin{split}
    \kappa^1 &= 2(\tfrac12 J_2A - \tfrac{\sqrt3}2 J_2B) \wedge (\tfrac12
    J_3A
    + \tfrac{\sqrt3}2 J_3B) - J_1A\wedge A \\
    &\qquad + (\tfrac12 J_2A - \tfrac{\sqrt3}2 J_2B) \wedge J_3A
    - (\tfrac12 J_3A + \tfrac{\sqrt3}2 J_3B) \wedge J_2A \\
    &= A\wedge J_1A + \tfrac32 J_2A \wedge J_3A\\
    &\qquad- \tfrac 32 J_2B \wedge J_3B
    + \sqrt3 J_2A \wedge J_3B + \sqrt3 J_3A \wedge J_2B.
  \end{split}
\end{equation*}
Giving
\begin{equation*}
  J_2\kappa^1 = - J_2A\wedge J_3A - \tfrac32 A \wedge J_1A
  + \tfrac 32 B \wedge J_1B - \sqrt3 A \wedge J_1B - \sqrt 3 J_1A \wedge
  B,
\end{equation*}
we then obtain
\begin{equation*}
  \begin{split}
    \kappa^1-J_2\kappa^1 &= \tfrac52 A\wedge J_1A + \tfrac52 J_2A \wedge
    J_3A -\tfrac32 B \wedge J_1B - \tfrac 32 J_2B \wedge J_3B \\
    &\qquad + \sqrt3(J_2A \wedge J_3B + J_3A \wedge J_2B + A \wedge J_1B +
    J_1A \wedge B),
  \end{split}
\end{equation*}
which is not proportional to
\begin{equation*}
  \omega_1 =  A\wedge J_1A + J_2A \wedge J_3A
  + B \wedge J_1B + J_2B \wedge J_3B.
\end{equation*}
In conclusion, $\QK{1+2}$ structures do not exist and the proof of
Theorem~\ref{thm:linear} is complete.

\section{Homogeneous descriptions of quaternionic hyperbolic space}
\label{sec:typ}

Being a non-compact symmetric space, the quaternionic hyperbolic
space~$\HH(n)$ admits a transitive (isometric) action of a solvable Lie
group, which is a proper subgroup of the full isometry group.  We thus see
that $\HH(n)$ has at least two homogeneous descriptions.  In this section
we study different homogeneous descriptions of~$\HH(n)$, listing all the
possible ones and finding their homogeneous types in some cases.  This is
preparation for~\S\ref{sec:qk3}, where we explicitly realise the
homogeneous structures of type~$\QK3$ on~$\HH(n)$.  A different model of
this construction is provided in~\S\ref{sec:example}.

\subsection{Transitive actions}

As a symmetric space, we have $\HH(n) = \SP(n,1)/(\SP(n)\times\SP(1)) =
G/K$.  A group $H$~acts transitively on $\HH(n)$ only if $H\backslash G/K$
is a point.  As $K$~is compact this implies that $H$ is a non-discrete
co-compact subgroup of the semi-simple group~$G$.  Such subgroups were
classified by Witte~\cite{Witte:cocompact} (cf.\ Goto \&
Wang~\cite{Goto-W:uniform}).

One begins by determining the standard parabolic subalgebras of $\lie g =
\sP(n,1)$ (cf.\ \cite[pp.~190--192]{Gorbatsevich-OV:structure}).  To be
concrete, we take $\sP(n,1)$ to be the set of quaternionic matrices that
are `anti-Hermitian' with respect to the bilinear form $B=\diag(\Id_{n-1},
\begin{smallmatrix} 0&1\\1&0 \end{smallmatrix})$, where $\Id_{n-1}$ is the
$(n-1)\times(n-1)$ identity matrix, thus
\begin{equation*}
  \sP(n,1) = \left\{\,
    \begin{pmatrix}
      \alpha & v_1 & v_2\\[0.2ex]
      -\bar v_2^T & a & b\\[0.2ex]
      -\bar v_1^T & c & -\bar a
    \end{pmatrix}
    :
    \begin{aligned}
      &\alpha\in\sP(n-1),\\
      &v_i\in\mathbb H^{n-1},\,a\in\mathbb H,\\
      &b,c\in\im\mathbb H
    \end{aligned}
    \,\right\}.
\end{equation*}
The maximal compact subgroup~$K$ has Lie algebra $\lie k =
\sP(n)\oplus\sP(1)$ with
\begin{gather*}
  \sP(n) = \left\{\,
    \begin{pmatrix}
      \alpha&v&v\\
      -\bar v^T&\beta&\beta\\
      -\bar v^T&\beta&\beta
    \end{pmatrix}
    :
    \begin{aligned}
      &\alpha\in\sP(n-1),\\
      &v\in\mathbb H^{n-1},\\
      &\beta\in\im\mathbb H
    \end{aligned}
      \,\right\},\\
  \sP(1) = \left\{\,
    \begin{pmatrix}
      0&0&0\\
      0&a&-a\\
      0&-a&a
    \end{pmatrix}
    : a\in\im\mathbb H \,\right\}.
\end{gather*}
Up to conjugation, $\sP(n,1)$~contains a unique maximal $\mathbb
R$-diagonalisable subalgebra $\lie a=\Span\{A\}$, with
$A=\diag(0,\dots,0,1,-1)$.  The set of roots $\Sigma$ corresponding
to~$\lie a$ is $\Sigma=\{\pm \lambda,\pm 2\lambda\}$, where $\lambda(A)=1$.
The set $\Theta=\{\lambda\}$ is a system of simple roots and the
corresponding positive roots system is~$\Sigma^+=\{\lambda,2\lambda\}$.

The general theory for $\lie g$ non-compact semi-simple says that
standard parabolic subalgebras correspond to subsets $\Psi$ of the
system of simple roots~$\Theta$ for a maximal $\mathbb
R$-diagonalisable subalgebra, as follows.  Let $[\Psi]$ be the
subset of~$\Sigma$ consisting of linear combinations of elements
of~$\Psi$.  Then a standard parabolic subalgebra $\lie
p(\Psi)=\lie l(\Psi) + \lie n(\Psi)$ of $\mathfrak{g}$ is
defined by
\begin{equation*}
  \lie l(\Psi) = \lie g_0 + \sum_{\mu\in [\Psi]} \lie g_\mu,\qquad
  \lie n(\Psi) = \sum_{\mu\in \Sigma^+\setminus [\Psi]} \lie g_\mu,
\end{equation*}
and each parabolic subalgebra of~$\lie g$ is conjugate to some $\lie
p(\Psi)$~\cite{Borel-T:reductive}.  The subalgebra $\lie n(\Psi)$ is
nilpotent, whilst $\lie l(\Psi)$ is reductive.  The latter may now be
decomposed as $\lie l(\Psi) = \lie l + \lie e + \lie a$, with $\lie
l$~semi-simple with all factors of non-compact type, $\lie e$~compact
reductive, and $\lie a$ the non-compact part of the centre of~$\lie
l(\Psi)$.  The decompositions
\begin{equation*}
  P(\Psi)^0 = LEAN,\qquad \lie p(\Psi) = \lie l + \lie e + \lie a
  + \lie n(\Psi)
\end{equation*}
are referred to as the \emph{refined Langlands decomposition} of the
parabolic subgroup $P(\Psi)$ and its Lie algebra in
Witte~\cite{Witte:cocompact} (cf.~\cite{Warner:harmonic1}).  Witte proved
the following theorem.

\begin{theorem}[Witte~\cite{Witte:cocompact}]
  \label{thm:Witte}
  Let $X$ be a normal subgroup of~$L$ and $Y$~a connected subgroup of~$EA$.
  Then there is a co-compact subgroup~$H$ of~$P(\Psi)$ whose component of
  the identity is $H^0 = XYN$.  Moreover, every non-discrete co-compact
  subgroup of~$G$ arises in this way.
\end{theorem}

For $G=\SP(n,1)$ we only have two choices for $\Psi$ and the corresponding
parabolic subalgebras have the following refined Langlands decompositions:
\begin{gather*}
  \lie p(\Theta) = \sP(n,1) + \{0\} + \{0\} + \{0\},\\
  \noalign{\smallskip} \lie p(\varnothing) = \{0\} + (\sP(n-1)+
  \sP(1)) + \lie a + (\lie n_1 + \lie n_2),
\end{gather*}
where
\begin{gather*}
  \sP(n-1)+\sP(1) = \left\{
    \begin{pmatrix}
      \alpha&0&0\\
      0&a&0\\
      0&0&a
    \end{pmatrix}
    : \alpha \in \sP(n-1),\, a \in \im\mathbb H
  \right\},\\
  \lie n_1 = \left\{
    \begin{pmatrix}
      0&0&v\\
      -\bar v^T&0&0\\
      0&0&0
    \end{pmatrix}
    : v \in \mathbb H^{n-1}
  \right\},\
  \lie n_2 =
  \left\{
    \begin{pmatrix}
      0&0&0\\
      0&0&b\\
      0&0&0
    \end{pmatrix}
    : b\in\im\mathbb H
  \right\},
\end{gather*}
the last being the $+1$- and $+2$-eigenspaces of~$\ad A$.  Note
that writing $N$ for the connected subgroup of~$\SP(n,1)$ with Lie
algebra $\lie n=\lie n_1 + \lie n_2$, we have that
$\SP(n,1)=KAN$ is the Iwasawa decomposition.

For the first case, Theorem~\ref{thm:Witte} says that for a co-compact~$H$,
the component of the identity~$H^0$ is either all of $\sP(n,1)$ or it is
trivial.  Thus the only transitive action coming from $\Psi=\varnothing$ is
that of the full isometry group $\SP(n,1)$ on~$\HH(n)$.

In the second case, there is much more freedom.  Each time we take a
connected subgroup~$Y$ of~$\SP(n-1)\SP(1)\mathbb R$ we get a corresponding
co-compact subgroup.  In order to get a transitive action on~$G/K$ note
that $\SP(n-1)\SP(1)$ is a subgroup of~$K$, so it is sufficient that the
projection $Y\to\SP(n-1)\SP(1)\mathbb R\to\mathbb R$ be surjective.
According to Theorem~\ref{thm:Witte}, $H^0$~is then~$YN$.  We thus have the
next result.

\begin{theorem}
  \label{thm:transitive}
  The connected groups acting transitively on~$\HH(n)$ are the full
  isometry group $\SP(n,1)$ and the groups $H=YN$, where $N$ is the
  nilpotent factor in the Iwasawa decomposition of $\SP(n,1)$ and $Y$~is a
  connected subgroup of $\SP(n-1)\SP(1)\mathbb R$ with non-trivial
  projection to~$\mathbb R$.
  \qed
\end{theorem}

The simplest choice is $Y=A$, this is then the description of~$\HH(n)$ as
the solvable group~$AN$.  One may determine a homogeneous type for this
solvable description as follows.

First we determine a quaternionic K\"ahler metric on~$AN$.  Taking the
decomposition $\lie a + \lie n_1 + \lie n_2$, a natural choice is
\begin{equation*}
  g(A,A) = \mu,\quad g(X_a,X_b) = \nu\delta_{ab},\quad g(V,V) = \norm v^2,
\end{equation*}
where $\norm{\cdot}$ is the Euclidean norm on $\mathbb H^{n-1}$ and
$\mu,\nu$ are positive constants.  A choice of a quaternionic
structure is then given by
\begin{equation*}
  J_1A = \kappa X_1, \quad J_1X_2 = X_3, \quad
  J_1\rho_1(v) = \rho_1(-vi),
\end{equation*}
where $\rho_1(v)=\left(
  \begin{smallmatrix}
    0&0&v\\ -\bar v^T&0&0\\ 0&0&0
  \end{smallmatrix}
\right)$ is the element of $\lie n_1$ corresponding to
$v\in\mathbb H^{n-1}$.  The compatibility condition
$g(A,A)=g(J_1A,J_1A)$ forces $\mu = \kappa ^2\nu$.  Taking
$v_1,\dots,v_{n-1}$ to be an orthonormal quaternionic basis of
$\mathbb H^{n-1}$, writing $V_i=\rho_1(v_i)$ and using
corresponding lower case letters to denote the left-invariant
basis dual to $A,X_1,\dots,V_1,J_1V_1,\dots,J_3V_{n-1}$, we have
\begin{equation*}
  \omega_1 = g(\cdot,J_1\cdot) = -\kappa\nu a\wedge x_1 - \nu x_2\wedge   x_3
  - \sum_{r=1}^{n-1} \left(v_r\wedge J_1v_r + J_2v_r\wedge J_3v_r\right),
\end{equation*}
etc.  The condition that the structure be quaternionic K\"ahler now reduces
to the requirement that $d\omega_1$ be a linear combination of $\omega_2$
and~$\omega_3$.  One computes
\begin{gather*}
  da = 0,\quad
  dv_r = -a \wedge v_r,\\
  dx_1 = -2a\wedge x_1 + 2\sum_{r=1}^{n-1} \left(v_r\wedge
    J_1v_r+J_2v_r\wedge  J_3v_r\right),
\end{gather*}
by using the fact that the exterior derivatives of the left-invariant
one-forms above are given by
\begin{equation*}
  d\alpha(B^*,C^*) = -\alpha([B,C]^*),
\end{equation*}
where $B^*$ is the vector field with one-parameter group $g\mapsto
\exp(tB)g$, $g\in G$, $t\in \mathbb R$, so $[B^*,C^*]=-[B,C]^*$.
Considering~$d\omega_1$ one finds that the structure is
quaternionic K\"ahler if and only if $\kappa\nu = -1$ and $\mu =
-\kappa = 1/\nu$.

The Levi-Civita connection~$\LC$ is given \cite[p.~183]{Besse:Einstein} by
\begin{equation*}
  2g(\LC_{B^*}C^*,D^*) =
  -\bigl\{g([B,C]^*,D^*)+g(B^*,[C,D]^*)+g(C^*,[B,D]^*)\bigr\}.
\end{equation*}
To find the homogeneous tensor $S = \LC - \Nt$, we need the canonical
connection~$\Nt$ for~$AN$.  But the latter is uniquely determined
\cite[p.~20]{Tricerri-V:homogeneous} by its value at $o \in \HH(n)$, where
we have
\begin{equation*}
  \Nt_{B^*}C^* = -[B,C]^*.
\end{equation*}
Note that $\Nt$ is the connection for which every left-invariant tensor
on~$AN$ is parallel \cite[p.~192]{Kobayashi-N:differential2}.  Working
at~$o$ and writing~$B$ for $B^*_o$, etc., we now have
\begin{equation}
  \label{eq:S1}
  2g(S_BC,D) = g([B,C],D) - g(B,[C,D]) - g(C,[B,D]).
\end{equation}
Note, for example, that this is skew-symmetric in $C$ and~$D$, confirming
that $S.g=0$.  To compute~$S$ explicitly, we first determine the Lie
brackets and find
\begin{align*}
  &g([B,C],D)=g(B,A_0)g(C,D)-g(C,A_0)g(B,D)\\
  &+\sum_{a=1}^32g(J_aB,C)g(D,J_aA_0)
  +\mu g(C,A_0)g(B,J_aA_0)g(D,J_aA_0)\\
  &\qquad\qquad -\mu g(B,A_0)g(C,J_aA_0)g(D,J_aA_0,D)\\
  &-2\mu \sum_{\sigma \in S_3}\varepsilon(\sigma) g(B,J_{\sigma (1)}A_0)
  g(C,J_{\sigma (2)}A_0) g(D,J_{\sigma(3)}A_0),
\end{align*}
where $A_0=A/\mu$.  Putting this into~\eqref{eq:S1} gives
\begin{equation}
  \label{gsbcd}
  \begin{split}
    &g(S_BC ,D)=-\sum_{a=1}^3g(B,J_aA_0)g(J_aC,D) \\
    &+g(D,A_0)g(C,B)-g(C,A_0)g(B,D)\\
    &\qquad+\sum_{a=1}^3\bigl(g(J_aB,C)g(D,J_aA_0)
    -g(D,J_aB)g(C,J_aA_0)\bigr)\\
    &{+}\mu\!\sum _{a=1}^3\bigl(g(C,A_0)g(B,J_aA_0)g(D,J_aA_0)
    {-}g(D,A_0)g(B,J_aA_0)g(C,J_aA_0)\bigr) \\
    &-\mu \sum_{\sigma \in S_3}\varepsilon(\sigma) g(B,J_{\sigma (1)}A_0)
    g(C,J_{\sigma (2)}A_0)g(D,J_{\sigma (3)}A_0).
  \end{split}
\end{equation}
The first line is in $\QK1$, the next two lines are a tensor in~$\QK3$.
With more work one finds that the last two lines are in $\QK{3+4}$ with
non-zero projection to each factor.  When $n>1$, the total contribution to
$\QK3$ is non-zero.  Although this computation is performed at~$o$, it
applies at each point of~$AN=\HH(n)$ and we have proved the next result.

\begin{proposition}
  For any fixed $\mu \in \mathbb R^+$, the quaternionic hyperbolic space $\HH(n)$
  admits a homogeneous quaternionic K\"ahler structure~$S$ of type
  $\QK{1+3+4}$ given by~\eqref{gsbcd}, corresponding to a description of
  $\HH(n)$ as a solvable group.
  \qed
\end{proposition}

We conclude that the class of a generic homogeneous structure on $\HH(n)$
is considerably larger than~$\QK3$.  Even without computation, one can see
that the solvable model necessarily has a non-trivial component
in~$\QK{1+2}$: the almost complex structures~$J_i$ are parallel with
respect to the canonical connection, but not with respect to the
Levi-Civita connection (since the scalar curvature is non-zero), so the
$\pi^a$ in equation~\eqref{spospnt} are not all zero.
Proposition~\ref{prop:QK3} then implies that there must also be a component
in~$\QK{4+5}$.

\subsection{Structures of type $\QK3$}
\label{sec:qk3}

We now determine the non-vanishing homogeneous structures of type~$\QK3$
on~$\HH(n)$.  Such a structure is given by a tensor
\begin{equation}
  \label{eq:QK3}
  S_XY = g(X,Y) \xi - g(\xi,Y) X +
  \sum_{a=1}^3 \bigl(g(\xi,J_aY) J_aX - g(X,J_aY) J_a\xi\bigr),
\end{equation}
where $\xi\ne0$~is a vector field satisfying $\Nt\xi=0$, i.e., $\LC\xi =
S\xi$.

We first consider the curvature term $R^S$ given by~\eqref{eq:RS}.
Using
\begin{gather*}
  S_X\xi = g(X,\xi)\xi - \norm\xi^2 X - \sum_{a=1}^3
  g(X,J_a\xi)J_a\xi,\qquad   S_\xi W = 0,\\
  \begin{split}
    S_XJ_1\xi  &= g(X,J_1\xi)\xi - \norm \xi^2 J_1X + g(X,\xi)J_1\xi \\
    &\qquad+ g(X,J_3\xi)J_2\xi - g(X,J_2\xi)J_3\xi,
  \end{split}\\
  S_{J_2\xi}W = 2\bigl\{ g(J_2\xi,W) \xi + g(J_3\xi,W) J_1\xi - g(\xi,W)
  J_2\xi - g(J_1\xi,W) J_3\xi \bigr\},
\end{gather*}
etc., we first expand
\begin{gather*}
  \begin{split}
    S_X(S_YW)
    &= g(Y,W) S_X\xi - g(\xi,W) S_XY\\
    &\qquad+ \sum_{a=1}^3\bigl(g(\xi,J_aW)S_X(J_aY) -
    g(Y,J_aW)S_X(J_a\xi)\bigr),
  \end{split}\\
  \begin{split}
    S_{S_XY}W &= g(X,Y)S_\xi W - g(\xi,Y)S_XW\\
    &\qquad + \sum_{a=1}^3 (g(\xi,J_aY)S_{J_aX}W - g(X,J_aY)S_{J_a\xi}W)
  \end{split}
\end{gather*}
and then anti-symmetrise to get
\begin{equation*}
  \begin{split}
    R^S_{XY}W
    &= S_Y(S_XW)-S_{S_YX}W -( S_X(S_YW)-S_{S_XY}W) \\
    &= - \norm\xi^2 \bigl\{ g(X,W)Y- g(Y,W)X\\
    &\qquad\qquad + \sum_{a=1}^3 \bigl(-g(X,J_aW)J_aY + g(Y,J_aW)J_aX\bigr)
      \bigr\}\\
    &\quad - 2\Cyclic_{123} \bigl\{ g(X,J_2Y)g(J_2W,\xi)\xi + g(X,J_2Y)
    g(W,\xi)J_2\xi \\
    &\qquad\qquad + g(X,J_3Y)g(J_1W,\xi)J_2\xi - g(X,J_1Y)g(J_3W,\xi)J_2\xi
      \bigr\}\\
    &\mathrel{=:} -\norm\xi^2 R^1_{XY}W + 2R^2_{XY}W.
  \end{split}
\end{equation*}
Note that $R^\ell_{XY}J_a = J_aR^\ell_{XY}$, for $\ell=1,2$ and
$a=1,2,3$.

From \eqref{eq:Rt-R} and the final lines of the proof of~\ref{prop:QK3}, we
have that
\begin{equation*}
  \Rt_{XY}W = \norm\xi^2 R^{\HH(n)}_{XY}W - R^S_{XY}W,
\end{equation*}
with
\begin{multline*}
  -R^{\HH(n)}_{XY}W = g(X,W) Y- g(Y,W) X \\
  + \sum_{a=1}^3\bigl(g(J_aX,W) J_aY - g(J_aY,W) J_aX + 2 g(J_aX,Y)
  J_aW\bigr) .
\end{multline*}
The first four terms here are exactly $R^1_{XY}W$, so
\begin{equation*}
  \Rt_{XY}W = - 2\norm\xi^2 \sum_{a=1}^3g(J_aX,Y) J_aW - 2R^2_{XY}W,
\end{equation*}
and we see in particular that $\Rt_{XY}\xi = 0$,
\begin{gather*}
  \Rt_{XY} J_1\xi = - 4\norm\xi^2\bigl\{ g(J_3X,Y)
  J_2\xi - g(J_2X,Y) J_3\xi\bigr\},\quad \text{etc.,}\\
  \text{and}\quad \Rt_{XY} Z = - 2\norm\xi^2 \sum _{a=1}^3g(J_aX,Y) J_aZ,
\end{gather*}
for $Z$ orthogonal to the quaternionic span of~$\xi$.  Thus $\Rt_{XY}$ acts
on $TM$ as an element of $\sP(1)$ in the representation $TM = \mathbb R
+ [S^2H] + [EH]$.  Also, for $Y$ orthogonal to~$\mathbb HX$, one
finds $\Rt_{XY}=0$, and for $\norm X^2=1/(2\norm\xi^2)$, we have that
\begin{equation*}
  \Rt_{XJ_aX}\xi=0,\quad \Rt_{XJ_aX}J_b\xi = -
  [J_a,J_b]\xi,\quad \Rt_{XJ_aX}Z = - J_aZ.
\end{equation*}
We write $\JJ_a$ for the element of $\sP(1)$ that acts as $J_a$ on the
factor~$[EH]$.

The corresponding homogeneous manifold $\widetilde G/\widetilde H$ has
(see~\S\ref{sec:AS})
\begin{equation*}
  \tilde{\lie h} = \sP(1),\quad \tilde{\lie g} = \tilde{\lie h} + T_oM,
\end{equation*}
with remaining Lie brackets (see~\eqref{corchete}) $[X,Y] = S_XY-S_YX +
\Rt_{XY}$, for $X,Y \in T_oM$.  On $T_oM$ we have
\begin{align}
  [Z_1,Z_2]& = 2\sum_{a=1}^3\bigl(g(J_aZ_1 ,Z_2)J_a\xi
  - \norm\xi^2g(J_aZ_1 ,Z_2)\JJ_a\bigr),\label{eq:ZZ}\\
  [\xi,Z] & = \norm\xi^2Z,\label{eq:xiZ}\\
  [\xi,J_a\xi] & = 2\norm\xi^2 J_a\xi - 2\norm\xi^4\JJ_a,\\
  [J_a\xi,Z] & = \norm\xi^2 J_aZ,\label{eq:JxiZ}\\
  [J_1\xi,J_2\xi] & = 4\norm\xi^2 J_3\xi - 2\norm\xi^4 \JJ_3, \quad
  \label{eq:IxiJxi}
\end{align}
for $Z,Z_1,Z_2$ orthogonal to $\mathbb H\xi$.

By Theorem~\ref{thm:transitive}, we need to identify this Lie
algebra structure as that of a subgroup $YN$ of
$\SP(n-1)\SP(1)\mathbb RN$, where $Y$ has non-trivial projection
to~$\mathbb R$.  Our holonomy algebra~$\lie h$ is isomorphic
to~$\sP(1)$, so the Killing form is negative definite on this
algebra and consequently it lies in $\sP(n-1)\oplus\sP(1)$. Indeed
it must lie in a subalgebra $\lie k_1\oplus\lie
k_2\subset\sP(n-1)\oplus\sP(1)$, with $\lie k_\ell\cong\sP(1)$.
Let $V$ and $W$ be the standard two-dimensional representations of
$\lie k_1$ and $\lie k_2$.  Then $\lie a + \lie n$ decomposes
as
\begin{equation*}
  \mathbb R + \bigl(\sum_{\ell\geqslant 0} c_\ell S^\ell V
  \otimes W\bigr)  + S^2W
\end{equation*}
under the action of $\lie k_1 \oplus \lie k_2$.  As the action of
the holonomy algebra only has a trivial summand of dimension~$1$,
we conclude that the projection of $\lie h$ to $\lie k_2$ is
non-zero.  Fitting the remaining representation to $S^2H +
(n-1)H$, we find (cf.~\cite[p.~110]{Fulton-Harris:rep}) that the
projection to $\lie k_1$ is zero.  Thus the holonomy algebra $\lie
h$ may be identified with $\lie k_2=\sP(1)$.

Comparing with~\S\ref{sec:AS}, we find that our symmetry group $\tilde G$
has Lie algebra
\begin{equation*}
  \tilde{\lie g} = \tilde{\lie h} + \lie m = \sP(1) + \lie a +
  \lie n_1 + \lie n_2,
\end{equation*}
and that $\lie n_2 \cong [S^2H]$.  Thus for each real~$\lambda$, there is
an $\ad$-invariant complement
\begin{equation}
  \label{com}
  \lie m_\lambda = \lie a
  + \lie n_1 + \lie p_\lambda, \qquad \text{where} \quad
  \lie p_\lambda=\left\{
    \begin{pmatrix} 0&0&0\\
      0&\lambda x&x\\
      0&0&\lambda x
    \end{pmatrix}
    : x\in\im\mathbb H \right\}.
\end{equation}
The relations \eqref{eq:xiZ}--\eqref{eq:IxiJxi} show that $\xi\in\lie a$
and $J_1\xi,J_2\xi,J_3\xi\in\lie p_\lambda $.  The holonomy action also
identifies $(\mathbb H\xi)^\bot $ with $\lie n_1$ and ensures that the map
$\rho_1\colon \mathbb H^{n-1}\to \lie n_1$ is quaternionic.  By considering
equation~\eqref{eq:ZZ}, we find that $\rho_1$ is a homothety and that
\begin{equation*}
  J_1\xi =
  \begin{pmatrix}
    0&0&0\\
    0&\norm\xi^2 i&i\\
    0&0&\norm\xi^2 i
  \end{pmatrix}
  ,\qquad
  \JJ_1 =
  \begin{pmatrix}
    0&0&0\\
    0&i&0\\
    0&0&i
  \end{pmatrix},
\end{equation*}
etc., once one has enforced the algebra requirement $[\JJ_1,\JJ_2]=2\JJ_3$.
Note that $\lambda=\norm\xi^2$.  One may now check this is consistent with
the remaining relations, including~\eqref{eq:JxiZ}.

\begin{theorem}
  \label{thm:prosixone}
  The non-vanishing homogeneous quaternionic K\"ahler structures of class
  $\QK3$ on $\HH(n)$ are parameterised by an element $\lambda \in \mathbb R^+$,
  corresponding to the $\ad$-invariant subspace $\lie m_\lambda$ described
  in~\eqref{com}.
  \qed
\end{theorem}

As a consequence of the previous results in the paper we now have the
characterisation of $\HH(n)$ stated in Theorem~\ref{thm:main}.

\subsection{A realisation of $\QK3$-structures}
\label{sec:example}

We finally give a realisation in the open ball model of $\HH(n)$ of the
family of $\QK3$-structures found in Theorem~\ref{thm:prosixone}.  Given a
fixed $c \in (-\infty,0)$, write $\tilde c=-c/4$ and consider the open
ball with radius $\rho_c=\sqrt{1/\tilde c}$ in $\HH(n)$,
\begin{equation*}
  B^n = \{ (q^0, q^1, \dots, q^{n-1}) \in \HH(n) : 1 -
  \tilde c \sum_{r=0}^{n-1} q^r {\bar q}^r > 0 \},
\end{equation*}
equipped with Watanabe's metric of negative constant quaternionic curvature
$c$ (see \cite[p.~134]{Watanabe:Q}) given, for $q^r = x^r + iy^r +jz^r + kw^r$,
$r = 0,\dots,n-1$, and $\qrad = \sum_{r=0}^{n-1} q^r {\bar
q}^r$, by
\begin{align*}
  g
  &= \frac1{(1 - \tilde c\qrad)^2}
  \{ \bigl( 1 - \tilde c \qrad \bigr) (dx^r dx^r + dy^r dy^r + dz^r
  dz^r + dw^r dw^r) \\ 
  &\quad + \tilde c \bigl( A_{rs}dx^r dx^s + B_{rs}dx^r
  dy^s + C_{rs}dx^r dz^s + D_{rs}dx^r dw^s \\
  &\qquad - B_{rs}dy^r dx^s + A_{rs}dy^r dy^s + D_{rs}dy^r dz^s -
  C_{rs}dy^r dw^s \\ 
  &\qquad - C_{rs}dz^r dx^s - D_{rs}dz^r dy^s + A_{rs}dz^r dz^s +
  B_{rs}dz^r dw^s \\ 
  &\qquad - D_{rs}dw^r dx^s + C_{rs}dw^r dy^s - B_{rs}dw^r dz^s +
  A_{rs}dw^r dw^s ) \},
\end{align*}
where
\begin{align*}
  A_{rs} & = x^r x^s + y^r y^s + z^r z^s + w^r w^s, \quad
  B_{rs} = x^r y^s - y^r x^s + z^r w^s - w^r z^s, \\
  C_{rs} & = x^r z^s - y^r w^s - z^r x^s + w^r y^s, \quad
  D_{rs} = x^r w^s + y^r z^s - z^r y^s - w^r x^s.
\end{align*}
Let $\upsilon$ be an almost quaternionic structure on $B^n$ admitting the
standard local basis
\begin{align*}
  J_1 & = \sum_{r=0}^{n-1} \Bigl( - \frac\partial{\partial
  x^r} \otimes d y^r + \frac\partial{\partial y^r} \otimes d x^r +
  \frac\partial{\partial z^r} \otimes d w^r -
  \frac\partial{\partial w^r} \otimes d z^r \Bigr), \\
  \noalign{\smallskip} J_2 & = \sum_{r=0}^{n-1} \Bigl( -
  \frac\partial{\partial x^r} \otimes d z^r - \frac\partial{\partial
  y^r} \otimes d w^r + \frac\partial{\partial z^r} \otimes d x^r +
  \frac\partial{\partial w^r} \otimes d y^r \Bigr),  \\
  \noalign{\smallskip} J_3 & = \sum_{r=0}^{n-1} \Bigl( -
  \frac\partial{\partial x^r} \otimes d w^r + \frac\partial{\partial
  y^r} \otimes d z^r - \frac\partial{\partial z^r} \otimes d y^r +
  \frac\partial{\partial w^r} \otimes d x^r \Bigr).
\end{align*}
Then, $(B^n,g,\upsilon)$ is a quaternionic K\"ahler manifold.  In fact, as
some computations show, conditions \eqref{jujdjt}, \eqref{jhjh}, and
\eqref{naji} are satisfied.  The Riemannian manifold $(B^n, g)$ is
homogeneous, hence complete.  Moreover, since $(B^n, g, \upsilon)$ is
connected, simply-connected and complete, it is a model of negative
constant quaternionic curvature.  We further look for a
homogeneous quaternionic K\"ahler structure $S$ satisfying~\eqref{eq:QK3}.
Since $\Nt\xi = 0$, we must find a vector field $\xi$ on $B^n$ such that
\begin{equation}
  \label{xi}
  \LC_X \xi = g(X,\xi)\xi- g(\xi,\xi)X - \sum_{a=1}^3g(X, J_a \xi) J_a\xi.
\end{equation}

The K\"ahler case \cite{Gadea-MA-MM:complex} suggests us the following
procedure to obtain~$\xi$.  Start with the Siegel domain
\begin{equation*}
  D_+ = \{ (\chi^0, \chi^1, \dots, \chi^{n-1}) \in \HH(n) : \re\chi^0 -
  \sum_{r=1}^{n-1} |\chi^r|^2 > 0 \}.
\end{equation*}
The Cayley transform $\varphi\colon B^n\to D_+$ has inverse~$\varphi^{-1}$
given by
\begin{equation*}
  (\chi^0, \chi^1, \dots, \chi^{n-1}) \mapsto
  \rho_c \Bigl(
  \frac{\chi^0 - 1}{\chi^0 + 1},  \frac{2\chi^1}{\chi^0 + 1}, \dots,
  \frac{2\chi^{n-1}}{\chi^0 + 1}\Bigr).
\end{equation*}
Writing $a_0=\re\chi_0$ and $\xi_{D_+} = 2\tilde c (a^0 - \sum_{r=1}^{n-1}
\abs{\chi^r}^2)\partial/\partial a^0$, we take $\xi=\varphi^{-1}_*\xi_{D+}$,
which is given explicitly by
\begin{align*}
  \xi= & \frac{\sqrt{\tilde c}\,(1 -\tilde c\qrad)}{\abs{q^0 -
  \rho_c}^2} \Bigl\{ \bigl( (x^0 -\rho_c)^2 - (y^0)^2 - (z^0)^2 -
  (w^0)^2 ) \frac\partial{\partial x^0} \\
  &\quad + 2(x^0 -\rho_c) \Bigl( y^0 \frac\partial{\partial y^0} + z^0
  \frac\partial{\partial z^0} + w^0 \frac\partial{\partial w^0}
  \Bigr)\\
   &\quad - \sum_{s=1}^{n-1} \Bigl( \bigl( (\rho_c - x^0) x^s + y^0 y^s +
   z^0 z^s + w^0 w^s \bigr) \frac\partial{\partial x^s} \\
   &\qquad + \bigl( (\rho_c-x^0) y^s - y^0 x^s - z^s w^0 + w^s z^0 \bigr)
   \frac\partial{\partial y^s} \\
   &\qquad + \bigl( (\rho_c-x^0)z^s - z^0 x^s + y^s w^0 - y^0 w^s \bigr)
   \frac\partial{\partial z^s}\\
   &\qquad + \bigl( (\rho_c-x^0) w^s - w^0 x^s - y^s z^0 - y^0 z^s \bigr)
   \frac\partial{\partial w^s} \Bigr) \Bigr\},
\end{align*}
satisfies $g(\xi,\xi) = \tilde c$ and equation~\eqref{xi}, so provides a
homogeneous structure of type~$\QK3$, with $\lambda=\tilde c$ in
Theorem~\ref{thm:prosixone}.

\providecommand{\bysame}{\leavevmode\hbox to3em{\hrulefill}\thinspace}
\providecommand{\MR}{\relax\ifhmode\unskip\space\fi MR }
\providecommand{\MRhref}[2]{%
  \href{http://www.ams.org/mathscinet-getitem?mr=#1}{#2}
}
\providecommand{\href}[2]{#2}

\end{document}